\documentclass[12pt,reqno]{amsart}

\usepackage[T1]{fontenc}
\usepackage{amsmath, amssymb, amsthm, amsthm, amscd, amsfonts}
\usepackage{mathtools}
\usepackage[dvipsnames]{xcolor}
\usepackage[alphabetic,lite,backrefs]{amsrefs}
\usepackage[colorlinks=true,hyperindex,linkcolor=magenta,pagebackref=true,citecolor=cyan]{hyperref}
\usepackage{algorithm, algorithmic}
\usepackage{caption, subcaption}
\usepackage{comment}
\usepackage{standalone}

\usepackage{extarrows}
\usepackage{tikz}
\usepackage{tikz-cd}
\usepackage{calrsfs}
\DeclareMathAlphabet{\pazocal}{OMS}{zplm}{m}{n}

\usepackage{fullpage}
\usepackage[shortlabels]{enumitem}

\usepackage{newcent}
\usepackage{fouriernc}

\newtheorem{theorem}{Theorem}[section]
\newtheorem{lemma}[theorem]{Lemma}

\newtheorem{corollary}[theorem]{Corollary}
\newtheorem{proposition}[theorem]{Proposition}

\theoremstyle{definition}

\newtheorem{definition}[theorem]{Definition}

\newtheorem{example}[theorem]{Example}

\newtheorem{convention}[theorem]{Convention}

\newtheorem{remark}[theorem]{Remark}

\newtheorem{theoremalpha}{Theorem}
\newtheorem{corollaryalpha}[theoremalpha]{Corollary}

\newtheorem{manualtheoreminner}{Theorem}
\newenvironment{manualtheorem}[1]{%
  \IfBlankTF{#1}
    {\renewcommand{\themanualtheoreminner}{\unskip}}
    {\renewcommand\themanualtheoreminner{#1}}%
  \manualtheoreminner
}{\endmanualtheoreminner}

\newtheorem{manualcorinner}{Corollary}
\newenvironment{manualcorollary}[1]{%
  \IfBlankTF{#1}
    {\renewcommand{\themanualcorinner}{\unskip}}
    {\renewcommand\themanualcorinner{#1}}%
  \manualcorinner
}{\endmanualcorinner}

\newcommand{\nn}{\mathbf{n}}

\newcommand{\bC}{\mathbf{C}}

\newcommand{\KK}{\mathbb{K}}
\newcommand{\NN}{\mathbb{N}}
\newcommand{\PP}{\mathbb{P}}
\newcommand{\QQ}{\mathbb{Q}}
\newcommand{\RR}{\mathbb{R}}
\newcommand{\ZZ}{\mathbb{Z}}
\newcommand{\CC}{\mathbb{C}}

\renewcommand{\O}{\mathcal{O}}

\newcommand{\cH}{\mathcal{H}}
\newcommand{\cE}{\mathcal{E}}
\newcommand{\cF}{\mathcal{F}}
\newcommand{\cI}{\mathcal{I}}

\DeclareMathOperator{\Pic}{Pic}
\newcommand{\PicR}[1]{\operatorname{Pic}_\RR\left(#1\right)}
\DeclareMathOperator{\Nef}{Nef}
\DeclareMathOperator{\Neg}{neg}
\DeclareMathOperator{\Eff}{Eff}
\newcommand{\NefR}[1]{\operatorname{Nef}_\RR\left(#1\right)}
\newcommand{\AmpR}[1]{\operatorname{Amp}_\RR\left(#1\right)}
\newcommand{\NR}[1]{N^1_\RR\left(#1\right)}

\DeclareMathOperator{\Proj}{Proj}

\DeclareMathOperator{\reg}{reg}
\newcommand{\regR}[1]{\operatorname{reg}_\RR\left(#1\right)}

\DeclareMathOperator{\sesh}{\mathbb{S}}
\newcommand{\ggen}[1]{\operatorname{\mathbb{G}}\left(#1\right)}
\newcommand{\ggenR}[1]{\operatorname{\mathbb{G}}_\RR\left(#1\right)}

\DeclareMathOperator{\conv}{conv}

\DeclareMathOperator*{\lcl}{lc.lim}
\DeclareMathOperator*{\ucl}{uc.lim}

\usetikzlibrary{decorations.pathreplacing}

\DeclareMathOperator{\Int}{int}
\DeclareMathOperator{\reint}{relint}

\DeclareMathOperator{\Sym}{Sym}
\DeclareMathOperator{\rank}{rank}
\DeclareMathOperator{\Amp}{Amp}

\newcommand{\BlX}{\widetilde{X}}

\newif\ifshownotes \shownotestrue 

\numberwithin{equation}{section}

\theoremstyle{theorem}
\newtheorem*{thm:asym-ideal-sheaf-reg}{Theorem~\ref{thm:ideal-sheaf-region-convergence}}
\newtheorem*{thm:asym-vect-bundle-reg}{Theorem~\ref{thm:vector-bundle-region-convergence}}

\title{Seshadri Regions and the Asymptotic Shape \\ of Multigraded Regularity}

\author{Juliette Bruce}
\address{Department of Mathematics, Dartmouth College, Hanover, NH}
\email{\href{mailto:juliette.bruce@dartmouth.edu}{juliette.bruce@dartmouth.edu}}
\urladdr{\url{https://juliettebruce.github.io}}

\author{Lauren Cranton Heller}
\address{Department of Mathematics, University of Nebraska, Lincoln, NE}
\email{\href{mailto:lch@math.berkeley.edu}{lheller2@nebraska.edu}}
\urladdr{\url{https://lcrantonh.github.io}}

\author{Mahrud Sayrafi}
\address{Department of Mathematics \& Statistics, McMaster University, Ontario, Canada}
\email{\href{mailto:mahrud@mcmaster.ca}{mahrud@mcmaster.ca}}
\urladdr{\url{https://ms.mcmaster.ca/~sayrafim/}}

\author{Alexandra Seceleanu}
\address{Department of Mathematics, University of Nebraska, Lincoln, NE}
\email{\href{mailto:aseceleanu@unl.edu}{aseceleanu@unl.edu}}
\urladdr{\url{https://aseceleanu.github.io}}


\begin{document}

\begin{abstract}
We introduce the Seshadri region of a subvariety, a convex region packaging the classical Seshadri constants with respect to every line bundle simultaneously. We develop the theory of Seshadri regions as a measure of positivity along subvarieties and apply it to determine asymptotic Castelnuovo–Mumford regularity for ideal powers and symmetric powers on smooth projective toric varieties.
\end{abstract}

\maketitle


\section{Introduction}\label{sec:intro}

A beautiful thread running between algebraic geometry and commutative algebra concerns bounds on the complexity of the defining equations of subvarieties of projective space \cites{GLP83,EisenbudGoto84,Lazarsfeld87,BEL91,MP18}. As succinctly stated by Cutkosky, Ein, and Lazarsfeld,
\begin{center}
	``\textit{... much of this material is clarified, and parts rendered transparent, when \\
    viewed through the lenses of vanishing theorems and intersection theory.}'' \cite{CEL01}
\end{center}
Another thread connecting these two communities is generalization of the classical algebra--geometry correspondence on projective space to other varieties, most prominently toric varieties \cites{Cox95,Mustata02,MS04,HMP10}. This article seeks to entwine these threads. We demonstrate that for subvarieties of a smooth projective toric variety, the complexity of defining equations is asymptotically governed by a convex body which is visible through the ``lenses'' of vanishing theorems and intersection theory.

The shift from projective space to toric varieties corresponds to a move from algebra over the standard $\ZZ$-graded polynomial ring $\KK[x_0,\ldots,x_n]$ with $\deg(x_i)=1$ to multigraded algebra over a $\ZZ^r$-graded polynomial ring $\KK[x_0,\ldots,x_n]$ with $\deg(x_i)$ in $\ZZ^r$. Hence measures of complexity of defining equations (e.g., degree, Castelnuovo--Mumford regularity, etc.) are no longer integers, but elements in $\ZZ^r$. A geometric consequence of this is that positivity is no longer measured against a single fixed ample divisor, the hyperplane class on $\PP^n$, but instead against the entire nef cone of $X$.

To be concrete, let $X$ be a projective algebraic variety defined over a field $\KK$. Given an ideal sheaf $\cI \subset \O_X$, consider the blowup $\pi\colon\BlX \to X$ along $\cI$ with exceptional divisor $E$. We define the Seshadri region of $\cI$ to be
\[
\sesh\left(\cI\right) \coloneqq \overline{\left\{ L\in \NR{X} \;\middle| \; \pi^*(L)-E \text{ is ample on $\BlX$} \right\}}.
\]
where the closure is taken in the Euclidean topology on $\NR{X}\cong \RR^r$ (see part~(3) of Proposition~\ref{prop:sesh-topology}). The Seshadri region is closed and convex, defined by finitely many rational hyperplanes when $\BlX$ is a Mori dream space, and on $X=\PP^n$ reduces to the reciprocal of Demailly's Seshadri constant from \cite{Demailly92} and to the s-invariant from \cite{CEL01}.

Heuristically, $\sesh(\cI)$ is a region with strong cohomological vanishing properties: it contains classes whose pullbacks dominate the exceptional divisor, allowing theorems like Fujita vanishing to hold on $\BlX$.  We apply this invariant to describe the multigraded Castelnuovo--Mumford regularity of an ideal sheaf $\cI$ on a smooth projective toric variety $X$.  Regularity is roughly a subset $\reg(\cI)$ of $\Pic(X)\cong N^1(X)\cong\ZZ^r$ capturing the cohomological complexity of~$\cI$ \cite{MS04}. Write $\regR{-}$ for $\reg(-)+\NefR{X}$, where $+$ denotes Minkowski sum. Our main result identifies the asymptotic regularity of powers of $\cI$ with the Seshadri region of~$\cI$.

\begin{manualtheorem}{\ref{main:thm-ideals}}
	Let $X$ be a smooth projective toric variety. If $\cI \subset \O_X$ is an ideal sheaf then
	\[
	\lim_{p\to \infty} \frac{\regR{\cI^p}}{p} = \sesh(\cI).
	\]
	where the limit is in the sense of Painlev\'{e}--Kuratowski convergence.
\end{manualtheorem}

A somewhat striking consequence of this theorem is that the limit of $\frac{\regR{\cI^p}}{p}$ is significantly nicer than $\reg(\cI^p)$ itself for any fixed $p$. For example, $\regR{\cI^p}$ is only convex when it contains a unique minimal element, which is unusual. This contrasts with the limit.

\begin{manualcorollary}{\ref{main:cor-convex}}
	Let $X$ be a smooth projective toric variety. If $\cI \subset \O_X$ is an ideal sheaf then $\lim_{p\to \infty} \frac{\regR{\cI^p}}{p}$ is convex.  Furthermore, if the blowup of $X$ along $\cI$ is a Mori dream space then it is defined by finitely many rational hyperplanes.
\end{manualcorollary}
Note that if $\cI$ defines a torus invariant subvariety of $X$ then the blowup of $X$ along $\cI$ is again a toric variety and hence a Mori dream space.

Corollary~\ref{main:cor-convex} is evidence that the asymptotics of multigraded regularity are governed by convex geometry. This mirrors similar constructions for other algebraic and geometric invariants, for instance Newton–Okounkov bodies, which capture volumes, base loci, and jet separation for graded linear series \cites{LM09,KK12,Anderson13,Ito13}; and Boij–Söderberg theory, which realizes Betti tables inside rational polyhedral cones \cites{BS08,EFW11}.

Further, computing multigraded regularity tends to be quite difficult, and there is no known finite algorithm. (See \cites{BCHS21,BCHS22} for some pathologies that can arise in the multigraded setting.) This is true even when $\cI$ defines a torus-invariant subvariety. Nevertheless we prove that $\sesh(\cI)$, and hence the asymptotic regularity of $\cI^p$, has a simple combinatorial description when $\cI$ defines a torus invariant subvariety.

This theorem generalizes significant work on the regularity of powers of homogeneous ideals \cites{Swanson97, CHT99, Kodiyalam00} and ideal sheaves of subvarieties in $\PP^n$. The most direct inspiration for our theorem comes from Cutkosky, Ein, and Lazarsfeld \cite{CEL01} who proved that if $\cI \subset \O_{\PP^{n}}$ is an ideal sheaf then
\[
	\lim_{p\to \infty} \frac{\reg\left(\cI^{p}\right)}{p} = \frac{1}{\varepsilon(\cI)}
\]
where $\varepsilon(\cI)$ is the Seshadri constant of $\cI$. Theorem~\ref{main:thm-ideals} implies this result, and our methods are largely inspired by theirs.  The proof rests on two main ideas: (i) due to the close relationship between regularity and global generation $\frac{\regR{\cI^p}}{p}$ is contained in $\sesh(\cI)$ for all~$p$ and (ii) elements in $\sesh(\cI)$ are positive enough that we may invoke Fujita vanishing on the blowup $\BlX$ to deduce the vanishing of $H^i(X, \cI^{p}(D))$ for divisors $D$ near $pL$ for $L\in\sesh(\cI)$.

Lastly, the analogy between blowups and projective bundles allows us to prove a similar result describing the asymptotic regularity of symmetric powers of vector bundles.

\begin{manualtheorem}{\ref{main:thm-bundles}}
	Let $X$ be a smooth projective toric variety. If $\cE$ is a vector bundle on $X$ and $F$ a Cartier divisor with $\O_{\PP(\cE)}(1)=\O_{\PP(\cE)}(F)$ then
	\[
	\lim_{p\to \infty} \frac{\regR{\Sym^{p}(\cE)}}{p} = \overline{\left\{ L\in \NR{X} \; \middle| \; \pi^*(L)+F \text{ is ample on $\PP(\cE)$} \right\}}
	\]
	where $\pi\colon\PP(\cE)\to X$ is the bundle map and the limit is in the sense of Painlevé--Kuratowski.
\end{manualtheorem}

As in the case of blowups we call this the Seshadri region of $\cE$, denoted $\sesh(\cE)$, which has many of the same properties as $\sesh(\cI)$. The proof of Theorem~\ref{main:thm-bundles} largely mimics the proof of Theorem~\ref{main:thm-ideals}, and in fact is easier in one respect: $\pi_*\O_{\PP(\cE)}(p)=\Sym^p(\cE)$ is almost definitional for vector bundles, whereas the analogous result for ideal sheaves and blowups is not.

\subsection*{Outline}  \S\ref{sec:setup-background} establishes conventions and reviews toric geometry (\S\ref{sub-sec:setup-toric}). \S\ref{sec:seshadri-region} introduces the Seshadri region of an ideal sheaf and proves basic results. \S\ref{sec:main-thm-regions} states our main result, Theorem~\ref{main:thm-ideals}, characterizing asymptotic behavior of multigraded regularity for powers of ideal sheaves, with examples and background on Painlevé--Kuratowski convergence. \S\ref{sec:ucl-sesh-global} and \S\ref{sec:lcl-sesh-fujita} prove Theorem~\ref{main:thm-ideals} and Corollary~\ref{main:cor-convex}: the first connects Seshadri regions to global generation, the second to Fujita vanishing; both proofs conclude in \S\ref{sub-sec:proof-of-main}. \S\ref{sec:parallel-bundles} proves Theorem~\ref{main:thm-bundles} for symmetric powers of vector bundles. \S\ref{sec:examples} computes explicit asymptotics for a point on a Fano surface. The Appendix covers Painlevé--Kuratowski convergence properties.

\subsection*{Acknowledgments}
We thank Christopher Manon for thoughtful questions which led us to Theorem~\ref{main:thm-bundles}, and Christine Berkesch, Daniel Erman, and David Eisenbud for advice and feedback. The computer software \textit{Macaulay2} \cite{M2} was vital in shaping our conjectures. AI tools were used to edit this paper.

This project began at the Simons Laufer Mathematical Sciences Institute (formerly MSRI) in Berkeley and continued at the Fields Institute in Toronto. JB was in residence at SLMath for the 2020-2021 academic year and LCH and AS for Spring 2024 (NSF DMS-1928930 and Alfred P.~Sloan Foundation G-2021-16778). The authors were partially supported by the National Science Foundation: JB under Award Nos.~DMS-1440140 and NSF MSPRF DMS-2002239; LCH under Award Nos.~NSF MSPRF DMS--2503497, DMS--2401482, and DMS--1901848; and AS under Award No.~DMS--2401482. MS was partially supported by the Doctoral Dissertation Fellowship at the University of Minnesota.

\section{Background \& Notation}\label{sec:setup-background}

Throughout we work over an algebraically closed field $\KK$. We denote the natural numbers by $\NN=\{0,1,2,\ldots\}$. If $T$ is a subspace of a $\KK$-vector space $V$ we write $\KK\langle T \rangle$ for the subspace of $V$ spanned by $T$. Similarly, we write $\RR_{\geq0}\langle T\rangle$ for the cone spanned by $T$.

\subsection{Algebraic Conventions}

We generally follow the conventions established in \cite{Lazarsfeld04I}; in particular, an (algebraic) variety is an irreducible, reduced, separated scheme of finite type over $\KK$. If $X$ is a complete variety, we write $\Pic(X)$ for the Picard group of $X$ (i.e., line bundles modulo isomorphism),  $N^1(X)$ for the Néron--Severi group of $X$ (i.e., Cartier divisors modulo numerical equivalence), and  $\NR{X}$ for the Néron--Severi group of $\RR$-divisors on $X$ (i.e., Cartier $\RR$-divisors modulo numerical equivalence). Note there is an isomorphism $\NR{X}\cong N^1(X)\otimes_{\ZZ}\RR$.

\subsection{Toric Varieties}\label{sub-sec:setup-toric}

Our main results, Theorems \ref{main:thm-ideals} and~\ref{main:thm-bundles}, will concern smooth projective toric varieties, although Section~\ref{sec:seshadri-region}  is applicable and of interest for all complete varieties, and we adopt that level of generality there. When working with toric varieties we adopt the standard conventions in \cite{CLS11}. If $X_{\Sigma}$ is a toric variety associated to a fan $\Sigma\subseteq \RR^n$ with a cone of full dimension then $\Pic(X_{\Sigma})$ is torsion free \cite{CLS11}*{Proposition 4.2.5}. If we further assume that $\Sigma$ has convex support then numerical equivalence and isomorphism of line bundles coincide and $\Pic(X_{\Sigma})\cong N^1(X_{\Sigma})$ \cite{CLS11}*{Theorem 6.3.15}.

All of these conditions are satisfied if we assume that $X$ is a smooth projective toric variety, which we generally will do. In this case $\Pic(X)$ is a finitely generated torsion-free abelian group (see \cite{CLS11}*{Theorems 4.2.1 and 4.2.5}) so $\Pic(X)\cong\ZZ^r$ for some $r\in\ZZ_{\geq0}$ and $\PicR{X}\cong \RR^r$. The set of nef divisors $\NefR{X}$ is a full dimensional cone in $\PicR{X}$ with interior the set $\Amp(X)_\RR$ of ample divisors \cite{CLS11}*{Theorem 6.3.22}. Note that $\NefR{X}$ may not be simplicial.

\begin{example}\label{ex:products-of-PPn}
  Given a positive vector $n=(n_1,\ldots,n_r)$ the product of projective spaces $\PP^\nn\coloneqq\PP^{n_1}\times \cdots \times \PP^{n_r}$ is a smooth projective toric variety. The Picard group of $\PP^\nn$ is isomorphic to $\ZZ^r$: given $d=(d_1,\ldots,d_r)\in \ZZ^r$ let $\O_{\PP^\nn}(d)\coloneqq \phi^*_1\O_{\PP^{n_1}}(d_1)\otimes \cdots \otimes \phi^*_r\O_{\PP^{n_r}}(d_r)$ where $\phi_i$ is projection onto the $i$th factor of $\PP^\nn$. The nef cone of $\PP^\nn$ is generated by $\O_{\PP^\nn}(e_i)$ for $i=1,\ldots,r$,  where the~$e_i$ are standard basis vectors in $\ZZ^r$. Thus $\Nef(\PP^\nn)$ is identified with $\ZZ^{r}_{\geq0}$ in $\ZZ^r$.
\end{example}

\begin{example}\label{eq:hirzebruch}
  The Hirzebruch surface $\cH_{a}\coloneqq \PP\left(\O_{\PP^1}\oplus\O_{\PP^1}(a)\right)$ with $a\geq 0$ is a smooth projective toric variety. The Picard group of $\cH_{a}$ is isomorphic to $\ZZ\langle E\rangle \oplus \ZZ\langle F\rangle \cong \ZZ^2$ where $F$ is the fiber of the map to $\PP^1$ and $E$ is the exceptional curve. Note that the classes of $F$ and $E+aF$ generate the nef cone of $\cH_a$.  The fan associated to $\cH_{a}$ is the complete two-dimensional fan with the four rays $\rho_0=(1,0)$, $\rho_1=(0,1)$, $\rho_2=(-1,a)$, and $\rho_3=(0,-1)$, shown below in the case of $\cH_2$. For each ray $\rho_i$ there is a corresponding prime torus invariant divisor $D_{\rho_i}$. These divisors satisfy the relations $D_{\rho_0}\sim D_{\rho_2}\sim F$, $D_{\rho_1}\sim E$, and $D_{\rho_3}\sim E+aF$.
  
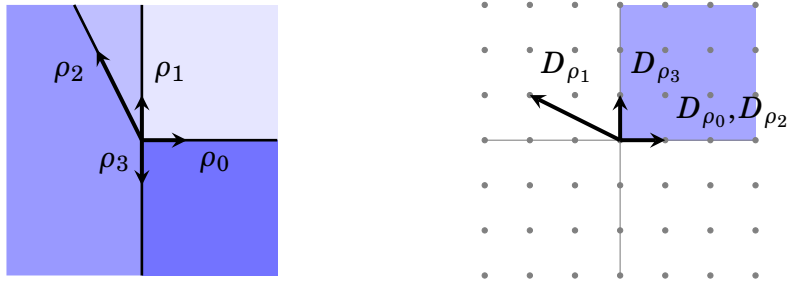
\begin{figure}[H]
  \begin{tikzpicture}[scale=.6]
    \path[use as bounding box] (-3,-3) rectangle (3,3);

    \fill[fill=blue!10] (0,0) -- (3, 0) -- ( 3,  3) -- ( 0, 3) -- cycle;
    \fill[fill=blue!25] (0,0) -- (0, 3) -- (-3/2,3) -- ( 0, 0);
    \fill[fill=blue!40] (0,0) -- (0,-3) -- (-3, -3) -- (-3, 3) -- (-3/2,3) -- (0,0);
    \fill[fill=blue!55] (0,0) -- (3, 0) -- ( 3, -3) -- ( 0,-3) -- cycle;

    \draw[line width=1pt,black] (0,0) -- (-3/2,3);
    \draw[line width=1pt,black] (0,0) -- ( 0, -3);
    \draw[line width=1pt,black] (0,0) -- ( 3,  0);
    \draw[line width=1pt,black] (0,0) -- ( 0,  3);

    \draw[line width=1.5pt,black,-stealth] (0,0) -- ( 1, 0) node[anchor=north west]{$\rho_0$}; 
    \draw[line width=1.5pt,black,-stealth] (0,0) -- ( 0,-1) node[anchor=south east]{$\rho_3$}; 
    \draw[line width=1.5pt,black,-stealth] (0,0) -- (-1, 2) node[anchor=north east]{$\rho_2$}; 
    \draw[line width=1.5pt,black,-stealth] (0,0) -- ( 0, 1) node[anchor=south west]{$\rho_1$}; 
  \end{tikzpicture}
  \hspace{1in}
  \begin{tikzpicture}[scale=.6]
    \path[use as bounding box] (-3,-3) rectangle (3,3);

    \fill[fill=blue!35] (0,0) -- (3,0) -- (3,3) -- (0,3) -- cycle;

    \draw [thin, gray,] (-3,0) -- (3,0); 
    \draw [thin, gray,] (0,-3) -- (0,3); 

    \foreach \x in {-3,...,3}{
      \foreach \y in {-3,...,3}{
        \node[circle,draw=gray, fill=gray,inner sep=0.7pt] at (\x,\y){ };
    }}

    \draw[line width=1.5pt,black,-stealth] (0,0) -- ( 1,0) node[anchor=south west]{$D_{\rho_0},D_{\rho_2}$};
    \draw[line width=1.5pt,black,-stealth] (0,0) -- (-2,1) node[anchor=south west]{$D_{\rho_1}$};
    \draw[line width=1.5pt,black,-stealth] (0,0) -- ( 0,1) node[anchor=south west]{$D_{\rho_3}$};
  \end{tikzpicture}
  \caption{Left: The fan of $\cH_2$. Right: $\Pic(\cH_2)$ with  $\Nef(\cH_2)$ (dark blue).}\label{fig:hirzebruch}
\end{figure}

\end{example}

We now describe one of the main topics of interest in this paper, multigraded regularity.  Given $\bC=(C_1,\ldots,C_t)$ a tuple of divisors and $s\in\ZZ^t$, write $s\cdot\bC$ for the sum $\sum_{i=1}^t s_iC_i$.

\begin{definition}\label{def:regularity}\cite{MS04}*{Definition~6.2}
  Let $X$ be a smooth projective toric variety. Fix a minimal generating set $\bC=(C_{1},\ldots,C_{t})$ for $\Nef(X)$ as a monoid. For $D\in \Pic(X)$, an $\O_{X}$-module $\cF$ is {\em $D$-regular}  if and only if for all $i>0$  the following vanishing occurs for $s\in\NN^t$:
  \[ H^i(X,\cF(D'))=0 \quad \text{for all} \quad D'\in \bigcup_{|s|=i}\left(D-s\cdot\bC+\Nef(X)\right). \]
  Collect all divisors $D$ satisfying this condition into the {\em multigraded regularity region}  \[\reg(\cF)\coloneqq\left\{D\in\Pic(X)\;\middle|\;\cF\text{ is $D$-regular}\right\}.\]
  We will eventually  take limits of regularity regions in $\PicR{X}$.  Since they are invariant under translation by $\Nef(X)$ by definition, it is reasonable to use this cone to create a closed region in $\PicR{X}$:
  \begin{equation}\label{regR}
  \regR{\cF}\coloneqq\reg(\cF)+\NefR{X}.
  \end{equation}
\end{definition}

\section{Seshadri Regions for Ideal Sheaves}\label{sec:seshadri-region}

In this section we introduce a multigraded substitute for the Seshadri constants introduced by Demailly and studied by Ein and Lazarsfeld over smooth projective surfaces \cites{Demailly92,EL92}.  In Section~\ref{sec:main-thm-regions} we will show that multigraded Seshadri regions describe the asymptotic behavior of multigraded regularity. Throughout this section we let $X$ be a complete algebraic variety.

\begin{definition}\label{def:sesh}
  Let $\cI \subseteq \O_X$ be an ideal sheaf and $\pi\colon\BlX \to X$ the blowup of $X$ along $\cI$ with exceptional divisor $E$. The \emph{Seshadri region of $\cI$ restricted to} a subspace $V \subset\NR{X}$ is
  \[
  \sesh_V\left(\cI\right) \coloneqq \left\{ L\in V \; \middle| \; \pi^*(L)-E \text{ is nef on $\BlX$} \right\} \subset\NR{X}.
  \]
  The \emph{Seshadri region} of $\cI$ is the Seshadri region of $\cI$ restricted to the $\RR$-span of $\NefR{X}$, and we denote this region by $\sesh(\cI)\coloneqq \sesh_{\NefR{X}}(\cI)$.
\end{definition}

The Seshadri region defined above is a globalization of the classical Seshadri constant. Classically, c.\,f.\,\cite{Lazarsfeld04I}*{Definition 5.1.1}, the Seshadri constant of the ideal sheaf $\cI$ of a closed point with respect to an ample class $H$ is defined to be the positive real number
  \[
  \varepsilon_H(\cI)=\max\left\{s\in \RR_{\geq 0} \; \middle| \;  \pi^*(H)-sE \text{ is nef }\right\}.
  \]

Setting $V$ in Definition \ref{def:sesh} to be the subspace $\RR\langle H\rangle$ spanned by $H$ makes $\sesh_H(\cI)$ a subset of $\RR$. The smallest element of this subset recovers the reciprocal of the Seshadri constant $\min \sesh_V(\cI) =\frac{1}{\varepsilon_H(\cI)}$. Equivalently this reciprocal is the distance from the origin to $\sesh(\cI)$ along the ray $\RR_{\geq0}\langle H\rangle$ (see Figure~\ref{fig:seshadri}).  Cutkosky, Ein, and Lazarsfeld refer to this reciprocal as the \emph{$s$-invariant} and define it for a wider class of ideal sheaves \cites{Paoletti94,Paoletti95,CEL01}.
\begin{center}
  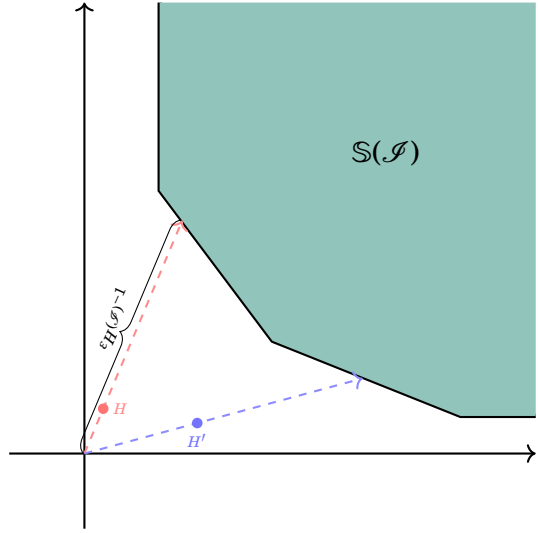
\begin{figure}[H]
  	\begin{tikzpicture}[scale=1]
  	  \draw[thick, ->] (-1, 0) -- (6, 0); 
	    \draw[thick, ->] (0, -1) -- (0, 6); 

	    \draw[ultra thick] (1,6)--(1,3.5) -- (2.5,1.5) -- (5,.5)--(6,.5);
	    \fill[PineGreen!40] (6,6)--(1,6)--(1,3.5) -- (2.5,1.5) -- (5,.5)--(6,.5)--(6,6);
	    \node at (4,4) {$\sesh(\cI)$};

	    \draw[thick, dashed, ->, red!45] (0,0) -- (1.3,3.1);
	    \draw[decorate,decoration={brace,amplitude=5pt,raise=0.5pt},yshift=0pt] (0,0) -- (1.3,3.1) node [midway,yshift=24pt, rotate = 66, xshift =-20pt]{\fontsize{5}{4}{ $\varepsilon_{H}(\cI)^{-1}$}};
	    \fill[red!55] (.25, 0.596155) circle (2pt);
	    \node[red!45, right] at (.25, 0.596155) {\fontsize{5}{4}{$ H$}};

	    \draw[thick, dashed, ->, blue!45] (0,0) -- (3.7,1);
	    \fill[blue!55] (1.5, 0.405405) circle (2pt);
	    \node[blue!55, below] at (1.5, 0.405405) {\fontsize{5}{4}{$H'$}};
    \end{tikzpicture}
	  \caption{The distance from the origin to the boundary of $\sesh(\cI)$ along a given ray is the reciprocal of the classical Seshadri constant $\varepsilon_{H}(\cI)$ for $H$ any ample class along the ray.}\label{fig:seshadri}
  \end{figure}
\end{center}

Thus the Seshadri region bundles all of the classical Seshadri constants. Since the reciprocal is more convenient for us to work with and since there is no danger of confusion between the classical Seshadri constant and the region described in Definition \ref{def:sesh}, we opted for the name Seshadri region.  We now establish its basic properties.

\begin{lemma}\label{lem:sesh-basic-properties}
  Let $X$ be a complete algebraic variety. Let $\cI \subseteq \O_X$ be an ideal sheaf and $\pi\colon\BlX \to X$ the blowup of $X$ along $\cI$ with exceptional divisor $E$. Fix a subspace $V\subset\NR{X}$.
  \begin{enumerate}
    \item\label{it:affine-transform} The linear map $\pi^*\colon\NR{X}\to\NR{\BlX}$ induces a bijection:
    \begin{equation}\label{eq:Sv-iso}
      \begin{tikzcd}[row sep = .5em, column sep = 3.5em]
        \sesh_V\left(\cI\right) \arrow[r, "\sim"] & \left(\pi^*(V)-E\right)\cap \NefR{\BlX}.
      \end{tikzcd}
    \end{equation}
    \item\label{it:translation} We have $\sesh_V(\cI) + (V\cap \NefR{X}) \subseteq \sesh_V(\cI)$.
    \item\label{it:inclusion} If $W \subseteq V \subset\NR{X}$ then $\sesh_W(\cI) \subseteq \sesh_V(\cI)$.
  \end{enumerate}
\end{lemma}
\begin{proof}
  (1) The map $\pi\colon\BlX\to X$ induces an $\RR$-linear inclusion $\pi^*\colon \NR{X}\to \NR{\BlX}$.  Its restriction given in \eqref{eq:Sv-iso} is surjective by Definition~\ref{def:sesh}.

  (2) Fix $L\in\sesh_V(\cI)$ and consider $L + L' \in L + (V\cap \NefR{X})$. Linearity of $\pi^*$ gives
  \[ \pi^*\big(L+L'\big)-E = \pi^*(L)-E + \pi^*(L'). \]
  Since $L\in \sesh_V(\cI)$ we know that $\pi^*(L)-E$ is nef. Since $\pi$ is proper and nef divisors pull back along proper maps  (see \cite{Lazarsfeld04I}*{Example 1.4.4(i)}),  $\pi^*(L')$ is nef on $\BlX$.  Thus $\pi^*(L+L')-E$ can be written as the sum of two nef divisors, so it is nef and $L+L'\in \sesh_V(\cI)$.

  (3) If $W\subseteq V$ then $\pi^*(W)-E\subseteq\pi^*(V)-E$ so $\sesh_W(\cI)\subseteq\sesh_V(\cI)$ by (1).
\end{proof}

An immediate corollary of Lemma~\ref{lem:sesh-basic-properties} is that the Seshadri region will be rational polyhedral region whenever the nef cone of the blowup $\BlX$ is a rational polyhedral cone. The question of when the nef cone is polyhedral has received significant attention \cites{HuKeel00,CT06,BCHM10,CT15,GK16,HKL18}.

\begin{corollary}\label{cor:sesh-poly}
  Let $X$ be a projective variety and $V \subset\NR{X}$ a rational subspace. If $\cI \subset \O_{X}$ is an ideal sheaf such that the blowup $\BlX$ is a Mori dream space then $\sesh_V(\cI)$ is a rational polyhedral region.
\end{corollary}

\begin{proof}
	If the blowup $\BlX$ of $X$ along $\cI$ is a Mori dream space then $\NefR{\BlX}$ is a rational polyhedral cone \cite{HuKeel00}. The claim follows from part~\eqref{it:affine-transform} of Lemma~\ref{lem:sesh-basic-properties} as the intersection of a rational affine subspace with a rational polyhedral region is itself rational polyhedral.
\end{proof}

When $X$ is a smooth projective toric variety with fan $\Sigma$, work of Di Rocco and others shows that one can compute Seshadri constants at torus fixed points (i.e., when $\cI$ is the ideal sheaf of such a point) combinatorially \cites{DiRocco99,BDHKKSS09}. We can provide a combinatorial description of the Seshadri region for arbitrary torus invariant subvarieties.

Write $D_0,\ldots, D_r$ for the torus-invariant divisors \emph{on the blowup $\BlX$} corresponding to the rays $\rho_0,\ldots,\rho_r$ in $\Sigma$, and $E$ for the exceptional divisor.  These divisors generate $N^1(\BlX)$.

\begin{proposition}\label{prop:sesh-toric}
  Let $X$ be a smooth projective toric variety.  Suppose that $\cI\subseteq \O_X$ is the ideal sheaf of the toric subvariety corresponding to a cone $\sigma\in\Sigma$ with $\dim\sigma>1$, and that $\sigma$ is spanned by rays $\rho_0,\ldots,\rho_j$.  Let $D_i$ be the strict transforms of the corresponding invariant divisors on $X$. Then the facets of $\sesh(\cI)$ correspond to (a subset of) the restrictions of the facets of $\NefR{\BlX}$ to the hyperplane of divisors $dE+\sum_{i=0}^rd_iD_i$ with $d=\sum_{i=0}^jd_i-1$.
\end{proposition}

\begin{proof}
 In this setting $\BlX$ is also a smooth projective toric variety and $N^1(\BlX)\cong\Pic(\BlX)$.  The pullback to $\BlX$ of a Cartier divisor on $X$ can be computed by evaluating its support function on the additional ray corresponding to $E$ (see \cite{CLS11}*{Proposition 6.2.7}).  For a prime torus invariant divisor on $X$ this gives either $D_i+E$ when $\rho_i\in\sigma$ or $D_i$ when $\rho_i\notin\sigma$.  Thus the set of pullbacks of divisors on $X$ is
  \(\pi^*\left(\Pic(X)\right)=\left\{dE+\sum_{i=0}^r d_iD_i \; \middle| \; d=\sum_{i=0}^j d_i\right\}\).  By part~\eqref{it:affine-transform} of Lemma~\ref{lem:sesh-basic-properties} the region $\sesh(\cI)$ is an affine transformation of $\left(\pi^*\left(N^1_{\mathbb{R}}(X)\right)-E\right)\cap\NefR{\BlX}$, whose facets are a subset of the restrictions of the facets of $\NefR{\BlX}$ to the hyperplane $\pi^*\left(\Pic(X)\right)-E$.  This hyperplane is exactly divisors $dE+\sum_{i=0}^rd_iD_i$ with $d=\sum_{i=0}^jd_i-1$.
\end{proof}

More recently, progress has been made on computing the Seshadri constants of line bundles at arbitrary points on toric varieties with additional structure \cites{Ito14,BDHK22, DR25}. It would be interesting to extend the polyhedral description from Proposition~\ref{prop:sesh-toric} in this vein. For more general ideal sheaves $\cI$, part~\eqref{it:affine-transform} of Lemma~\ref{lem:sesh-basic-properties} makes it clear that $\sesh_{V}(\cI)$ may not be rational polyhedral.

\begin{example}
	If we take $\cI=\O_{X}$ then the blowup $\pi\colon\BlX \to X$ is just the identity map, so $\sesh_{V}(\O_X)=V\cap \NefR{X}$. Therefore $\sesh(\O_X)$ will not be rational polyhedral whenever $\NefR{X}$ is not. Examples include some abelian varieties and products of curves \cites{Bauer98,Rabindranath19}.
\end{example}

\begin{example}
	If $X=\PP^1\times\PP^1$ and $\cI \subset \O_{\PP^1\times \PP^1}$ is the ideal sheaf of 10 very general points then $\varepsilon_D(\cI) = \frac{1}{\sqrt{5}}$ where $D$ is the class of $\O_{\PP^1\times\PP^1}(1,1)$ \cite{LU26}*{Theorem~1}.  In particular, even when $\sesh(\cI)$ is defined by finitely many inequalities (here one), they need not exist over $\QQ$.
\end{example}

Despite these pathologies the Seshadri region $\sesh_V(\cI)$ is closed and convex. Its relative interior is nonempty when $V$ contains an ample divisor on $X$, in which case we can describe it in terms of ample divisors on $\BlX$. We use relative interior since $\sesh_{V}(\cI) \subset V$ and $V$ may be a proper subspace, in which case $\Int(\sesh_{V}(\cI))$ will be empty.  The relative interior of a subset $S \subset \RR^r$ is the interior of $S$ as a subset of its affine span.

\begin{proposition}\label{prop:sesh-topology}
  Assume $X$ is a projective variety. Let $\cI \subseteq \O_X$ be an ideal sheaf and let $\pi\colon\BlX \to X$ be the blowup of $X$ along $\cI$ with exceptional divisor $E$. Fix a subspace $V\subset\NR{X}$ and assume that $V\cap\Amp_\RR(X)\neq\varnothing$.
  \begin{enumerate}
    \item The Seshadri region of $\cI$ restricted to $V$ is a nonempty closed convex set with respect to the Euclidean topology on $\NR{X}$.
    \item Its relative interior is
    \[
    \reint\left(\sesh_{V}\left(\cI\right)\right) = \left\{L\in V \; \middle| \; \pi^*(L)-E \text{ is ample on $\BlX$} \right\} \subset\NR{X}.
    \]
    \item  When $V=\NefR{X}$ the interior of $\sesh(\cI)$ is nonempty and
    \[
    \sesh\left(\cI\right) = \overline{ \left\{L\in\PicR{X} \; \middle| \; \pi^*(L)-E \text{ is ample on $\BlX$} \right\} }.
    \]
  \end{enumerate}
\end{proposition}

\begin{proof}
	 Let $F\colon\NR{X} \to \NR{\BlX}$ be the affine linear map defined by $F(L)=\pi^*(L)-E$.  Then from part~\eqref{it:affine-transform} of Lemma~\ref{lem:sesh-basic-properties}
	 \[
   \sesh_V(\cI)= \left\{ L \in V \;\middle|\; \pi^*(L) - E \in \NefR{\BlX}\right\} = V \cap F^{-1}\left(\NefR{\BlX}\right).
	 \]
	 Since $X$ is projective the blowup $\BlX$ is also projective, so Kleiman’s criterion implies that $\NefR{\BlX}$ is a closed, convex cone of full dimension in $\NR{\BlX}$, with interior $\AmpR{\BlX}$ \cite{Lazarsfeld04I}*{Theorem~1.4.23}.

   As $F$ is an affine map $F^{-1}\left(\NefR{\BlX}\right)$ is closed, convex, and has nonempty interior.  Since $V$ is a subspace of $\RR^r$ the intersection $\sesh_V(\cI)$ is also closed and convex.  Since there exists $L\in V\cap\AmpR{X}$ and $-E$ is ample for $\pi$ \cite{stacks-project}*{Lemma 02NS}, the divisor $\pi^*(\ell L)-E$ is ample for some $\ell\in\NN$.  Thus $\sesh_V(\cI)$ is nonempty.

   Turning toward part~(2), the existence of $\ell L\in V\cap F^{-1}(\AmpR{\BlX})$ together with \cite{Rock70}*{Corollary 6.5.1 and Theorem 6.7} imply that
\begin{align*}
   V\cap F^{-1}\left(\AmpR{\BlX}\right) &= V \cap F^{-1}\left(\reint\left(\NefR{\BlX}\right)\right)\\
   &= V \cap \reint\left( F^{-1}\left(\NefR{\BlX}\right) \right)\\
   &= \reint\left(V \cap F^{-1}\left(\NefR{\BlX}\right) \right) = \reint\left(\sesh_{V}\left(\cI\right)\right)
\end{align*}
   For part (3) we must show that a point $L\in\sesh(\cI)$ is contained in the closure of the interior of $\sesh(\cI)$, since we already know $\sesh(\cI)$ is closed and nonempty.  By part \eqref{it:translation} of Lemma~\ref{lem:sesh-basic-properties} the translated full dimensional cone $L+\NefR{X}$ is contained in $\sesh(\cI)$.  Thus $L$ is contained in the closure of the open subset $L+\AmpR{X}$ of $\sesh(\cI)$.
\end{proof}

\section{The Main Theorem \& Convergence of Regions}\label{sec:main-thm-regions}

For the remainder of the paper we specialize to the case when $X$ is a smooth projective toric variety and fix the following notation.

\begin{convention}\label{notation:blowup}
	Let $X$ be a smooth projective toric variety of dimension $n$, and suppose $\bC=(C_{1},\ldots,C_{t})$ is an ordered tuple of divisors such that $C_{1},\ldots,C_{t}$ is a minimal generating set of $\Nef(X)$.  Let $\cI \subseteq \O_{X}$ be an ideal sheaf and $\pi\colon\BlX \to X$ the blowup of $X$ along $\cI$ with exceptional divisor $E$.
\end{convention}

We can now state the main theorem, connecting Seshadri regions and the asymptotics of multigraded regularity.  As mentioned in Section~\ref{sub-sec:setup-toric}, under the conditions of Convention~\ref{notation:blowup} we have $N^1(X)\cong \Pic(X)$, allowing us to compare $\sesh(\cI)$ and $\regR{\cI}$ in the same vector space.

\begin{theoremalpha}\label{main:thm-ideals}
  Let $X$ be a smooth projective toric variety. If $\cI \subset \O_X$ is an ideal sheaf then
	\[
	\lim_{p\to \infty} \frac{\regR{\cI^p}}{p} = \sesh(\cI).
	\]
	where the limit is in the sense of Painlevé--Kuratowski convergence by Definition~\ref{def:limits} below.
\end{theoremalpha}

\begin{corollaryalpha}\label{main:cor-convex}
	Let $X$ be a smooth projective toric variety. If $\cI \subset \O_X$ is an ideal sheaf then $\lim_{p\to \infty} \frac{\regR{\cI^p}}{p}$ is convex. Furthermore, if the blowup of $X$ along $\cI$ is a Mori dream space then the limit of scaled regularity regions is defined by finitely many rational hyperplanes.
\end{corollaryalpha}

\begin{remark}
  If $\overline{\cI}$ is the integral closure of an ideal sheaf $\cI \subset \O_{X}$ then $\reg(\cI)$ and $\reg(\overline{\cI})$ can be different. However Theorem~\ref{main:thm-ideals} implies that asymptotically the scaled regularities of $\cI^p$ and $\overline{\cI}^p$ agree, since the normalization of the blowup along $\overline{\cI}$ is the normalization of the blowup of $X$ along $\cI$, the pullbacks of the exceptional divisors to this common normalization agree, and proper surjective maps preserve nefness.

  Put differently, the asymptotic regularity of $\cI$ only depends on the integral closure of $\cI$, even if $\cI$ is not the ideal sheaf of a normal toric subvariety.
\end{remark}

For the remainder of this section we recall the notion of Painlev\'e--Kuratowski convergence for sequences of subsets in $\RR^r$ and give examples.  Let $d(x,y)$ denote the standard Euclidean metric throughout. For $A\subseteq \RR^r$ and $x\in \RR^r$ define
\[
d(x,A)\coloneqq \inf_{y\in A} d(x,y).
\]

\begin{definition}\label{def:limits}
	Let $\{A_{p}\}$ be a sequence of subsets of $\RR^r$.
	\begin{itemize}
		\item The \emph{lower closed limit} of $\{A_{p}\}$ is
		\[
		\lcl_{p\to\infty} A_{p} \coloneqq \left\{x\in \RR^r \;\middle|\; \limsup\limits_{p\to\infty}d(x,A_{p})=0\right\}.
		\]
		\item The \emph{upper closed} limit of $\{A_{p}\}$ is
		\[
		\ucl_{p\to\infty} A_{p} \coloneqq \left\{x\in \RR^r \;\middle|\; \liminf\limits_{p\to\infty}d(x,A_{p})=0\right\}.
		\]
		\item The limit of $\{A_{p}\}$ exists and equals $L\subseteq \RR^r$ if
		\[
		\lim_{p\to\infty}A_{p}\coloneqq \ucl_{p\to\infty} A_{p} =\lcl_{p\to\infty} A_{p}=L.
		\]
	\end{itemize}
If the limit exists we say that the sequence of subsets $\{A_{p}\}$ {\em converges to} $L$ in the sense of Painlev\'e--Kuratowski.
\end{definition}

When we feel there is relatively little room for confusion we will simply write $\lcl \{A_{p}\}$ for $\displaystyle \lcl_{p\to\infty} A_{p}$, $\ucl \{A_{p}\}$ for $\displaystyle \ucl_{p\to\infty}$, and $\lim \{A_{p}\}$ for $\displaystyle \lim_{p\to\infty}A_{p}$.  While standard in certain areas of mathematics, it seems that Painlev\'e--Kuratowski convergence has found less use in algebra and geometry.  We point an unfamiliar reader to Appendix \ref{appendix} for a number of basic properties and their proofs.

\begin{example}\label{ex:translating-cone}
	Fix a vector $v \in \RR^r_{\geq0}$. Consider the sequence of subsets $\{A_{p}\}$ of $\RR^r$ where $A_{p}=pv+\RR^r_{\geq0}$. Since $v\in \RR^r_{\geq0}$ this is a nested decreasing sequence $A_{1}\supseteq A_{2} \supseteq A_{3} \supseteq \cdots$, so Lemma~\ref{lem:region-chains} implies the limit of this sequence exists and is equal to $\bigcap \overline{A}_{p}$.  However, it is clear that this intersection is empty so long as $v\neq 0$.

  If we scale the above regions by letting $\widehat{A}_{p}=\frac{A_p}{p}$ we see that the limit is more interesting. In particular, $\widehat{A}_{p}=v+\RR^r_{\geq 0}$ for all $p$, so $\lim \{\widehat{A}_{p}\}=v+\RR^r_{\geq0}$ by part~(5) of Lemma~\ref{lem:region-limit-properties}. This can be generalized to showing that if $f(p)$ is a linear function of $p$ and $A_{p}=f(p)v+V$ for some closed pointed cone $V$ in $\RR^r_{\geq 0}$ then $\lim \{\widehat{A}_{p}\}=mv+V$ where $m$ is the slope of $f(p)$.
\end{example}

We illustrate Theorem \ref{main:thm-ideals} with an explicit example. Example \ref{ex:blowup-P1xP1} is deceptively simple: a rare instance where one is able to determine the regularity regions of all powers of $\cI$. Section~\ref{sec:examples} presents a more elaborate example which better reflects the strength of the result.

\begin{example}\label{ex:blowup-P1xP1}
  Let $X=\PP^1\times\PP^1$, let $I$ be the ideal of a torus-invariant point, and let $\cI=\widetilde I$.  After identifying $\Pic(X)$ with $\ZZ^2$ by $(a,b)\mapsto\O_{\PP^1}(a)\boxtimes\O_{\PP^1}(b)=\O_X(a,b)$ the nef cone of $X$ is $\RR_{\geq0}^2$, so we can take $\bC$ to be the standard basis $(e_1,e_2)$.  In these coordinates $\sesh(\cI)=(1,1)+\RR_{\geq 0}^2$.

  Since $I$ is a complete intersection, the resolution of $I^p$ is an Eagon--Northcott complex with twists
  \[0\to\bigoplus_{a+b=p-1}S(-a-1,-b-1)\to\bigoplus_{a+b=p}S(-a,-b)\to I^p \to 0\]
  for $a$ and $b$ nonnegative integers.    By \cite{MS04}*{Corollary~7.3}, and using their convention of a map $\phi\colon\{1,2,3\}\to\{1,2\}$, we have
  \begin{align*}
    \regR{I^p}&\supseteq \bigcup_{\phi\colon[3]\to[2]}\left(\left[\bigcap_{a+b=p}(a,b)+\RR_{\geq0}^2\right]\cap\left[\bigcap_{a+b=p-1}(a+1,b+1)-e_{\phi(1)}+\RR_{\geq0}^2\right]\right)\\
    &=\bigcup_{\phi\colon[3]\to[2]}\left(\left[(p,p)+\RR_{\geq0}^2\right]\cap\left[(p,p)-e_{\phi(1)}+\RR_{\geq0}^2\right]\right)\\
    &=\bigcup_{\phi\colon[3]\to[2]}\left((p,p)+\RR_{\geq0}^2\right)\\
    &=(p,p)+\RR_{\geq0}^2
  \end{align*}
  as $\regR{S}=\RR_{\geq0}^2$.  Conversely, if $(a,b)\notin(p,p)+\RR_{\geq0}^2$ then $I^p_{\geq(a,b)}$ has generators in either degree $(p,b)$, when $a<p$, or $(a,p)$, when $b<p$.  Therefore $I^p$ is not $(a,b)$-regular by the contrapositive of \cite{MS04}*{Theorem~5.4}, so we have established that $\regR{I^p}=(p,p)+\RR_{\geq0}^2$.

  Since $I$ is generated by a regular sequence, its powers are saturated \cite{ZariskiSamuel}*{Lemma 5, Appendix 6}, so $\cI^m=\widetilde{(I^m)}$ and $\regR{\cI^p}=\regR{I^p}=(p,p)+\RR_{\geq0}^2$. As in Example \ref{ex:translating-cone} we conclude that
  \[\textstyle
	\lim_{p\to \infty} \frac{\regR{\cI^p}}{p} = (1,1)+\RR_{\geq0}^2=\sesh(\cI).
  \]
\end{example}

Multigraded regularity and Seshadri regions are invariant under translation by certain pointed cones (see Definition~\ref{def:regularity} and part~\eqref{it:translation} of Lemma~\ref{lem:sesh-basic-properties}), so it will be relevant that this property is preserved in their limits.

\begin{lemma}\label{lem:limits-semi-groups}
	Let $V\subseteq \RR^r$ be a closed pointed cone whose vertex is the origin. Let $\{B_p\}$ be a sequence of subsets of $\RR^r$ such that $\lim\{B_p\}=L$. If $B_p\subseteq L$ for all $p$ and the Minkowski sum $L+V\coloneqq\left\{x+y \;\middle|\; x\in L, \; y\in V\right\}$ is closed then $\lim\{B_p+V\}$ exists and is equal to $L+V$.
\end{lemma}

Lemma~\ref{lem:limits-semi-groups} allows us to easily compute the limits of regions whose component-wise minimal elements converge.

\begin{example}\label{ex:region-limit-staircase}
  Fix nonnegative integers $a,b\in \ZZ_{>0}$.  Consider the sequence of subsets $\{A_{p}\}$ of $\RR^{2}$ where
  $ A_{p}\coloneqq \bigcup\limits_{i=0}^{p}\left[\left(a(p-i),bi\right)+\RR_{\geq0}^2\right]$.
  As in Example~\ref{ex:translating-cone} the limit of the sequence $\{A_{p}\}$ is uninteresting since it is empty, so we consider the scaled sequence $\{\widehat{A}_{p}\}$ with $\widehat{A}_{p}\coloneqq \tfrac{A_{p}}{p}$.  Defining \(B_p=\left\{(a(p-i),bi)\;\middle|\; i=0,1,\ldots,p\right\}\) we see that $\lim\{\hat B_p\}$ for $\hat B_p\coloneqq\tfrac{B_p}{p}$ is the convex hull $\conv\left((a,0),(0,b)\right)$. Thus $\lim\left\{\widehat{A}_p\right\}=\conv\left((a,0),(0,b)\right)+\RR^{2}_{\geq0}$ by Lemma~\ref{lem:limits-semi-groups}.
	
\begin{figure}[H]
\newcommand{\makegrid}{
  \path[use as bounding box] (-1.45,-1.25) rectangle (9.45,7.25);
  \draw[-,  semithick] (-1,0)--(9,0);
  \draw[-,  semithick] (0,-1)--(0,7);
  \node at (-1,3) (a) {\scriptsize $(0,3)$};
  \node at (7,-1) (a) {\scriptsize $(7,0)$};
  \fill[fill=Black] (6.95,-.15) rectangle (7.05,0.15);
  \fill[fill=Black] (-.15,2.95) rectangle (0.15,3.05);
    \draw[-,  dashed] (0,3)--(7,0);
}

\begin{tikzpicture}[scale=.4]
  \path[fill=PineGreen!45] (7, 0)--(7, 3)--(0, 3)--(0,7)--(9,7)--(9,0)--(7,0);
  \makegrid

  \draw[->, ultra thick,PineGreen] (7,0)--(9,0);
  \draw[-, cap=round,ultra thick,PineGreen] (7, 0)--(7, 3)--(0, 3);
  \draw[->, ultra thick,PineGreen] (0,3)--(0,7);

  \fill[TealBlue,fill=PineGreen] (7,0) circle (5pt);
  \fill[TealBlue,fill=PineGreen] (0,3) circle (5pt);
      \node at (5,4.55) (a) {\scriptsize $\widehat{A}_1$};
\end{tikzpicture}\quad
%
%
\begin{tikzpicture}[scale=.4]
  \path[fill=PineGreen!45] (7,0)--(7,1)--(4.66667,1)--(4.66667,2)--(2.33333,2)--(2.33333,3)--(0,3)--(0,7)--(9,7)--(9,0)--(7,0);
  \makegrid

  \draw[->, ultra thick,PineGreen] (7,0)--(9,0);
  \draw[-, cap=round,ultra thick,PineGreen] (7,0)--(7,1)--(4.66667,1)--(4.66667,2)--(2.33333,2)--(2.33333,3)--(0,3);
  \draw[->, ultra thick,PineGreen] (0,3)--(0,7);

  \fill[TealBlue,fill=PineGreen] (7,0) circle (5pt);
  \fill[TealBlue,fill=PineGreen] (0,3) circle (5pt);
  \fill[TealBlue,fill=PineGreen] (4.66667,1) circle (5pt);
  \fill[TealBlue,fill=PineGreen] (2.33333,2) circle (5pt);

      \node at (5,4.55) (a) {\scriptsize $\widehat{A}_3$};
\end{tikzpicture}\quad
\begin{tikzpicture}[scale=.4]
  \path[fill=PineGreen!45] (7,0)--(7,.6)--(5.6,.6)--(5.6,1.2)--(4.2,1.2)--(4.2,1.8)--(2.8,1.8)--(2.8
      ,2.4)--(1.4,2.4)--(1.4,3)--(0,3)--(0,7)--(9,7)--(9,0)--(7,0);
  \makegrid

  \draw[->, ultra thick,PineGreen] (7,0)--(9,0);
  \draw[-, cap=round,ultra thick,PineGreen] (7,0)--(7,.6)--(5.6,.6)--(5.6,1.2)--(4.2,1.2)--(4.2,1.8)--(2.8,1.8)--(2.8
      ,2.4)--(1.4,2.4)--(1.4,3)--(0,3);
  \draw[->, ultra thick,PineGreen] (0,3)--(0,7);

  \fill[TealBlue,fill=PineGreen] (7,0) circle (5pt);
  \fill[TealBlue,fill=PineGreen] (0,3) circle (5pt);
  \fill[TealBlue,fill=PineGreen] (5.6,.6) circle (5pt);
  \fill[TealBlue,fill=PineGreen] (4.2,1.2) circle (5pt);
  \fill[TealBlue,fill=PineGreen] (2.8,1.8) circle (5pt);
  \fill[TealBlue,fill=PineGreen] (1.4,2.4) circle (5pt);

      \node at (5,4.55) (a) {\scriptsize $\widehat{A}_5$};
\end{tikzpicture}\quad

\begin{tikzpicture}[scale=.4]
  \path[fill=PineGreen!45] (7,0)--(7,.3)--(6.3,.3)--(6.3,.6)--(5.6,.6)--(5.6,.9)--(4.9,.9)--(4.9,1.2
      )--(4.2,1.2)--(4.2,1.5)--(3.5,1.5)--(3.5,1.8)--(2.8,1.8)--(2.8,2.1)--(2.1
      ,2.1)--(2.1,2.4)--(1.4,2.4)--(1.4,2.7)--(.7,2.7)--(.7,3)--(0,3)--(0,7)--(9,7)--(9,0)--(7,0);
  \makegrid

  \draw[->, ultra thick,PineGreen] (7,0)--(9,0);
  \draw[-, cap=round,ultra thick,PineGreen] (7,0)--(7,.3)--(6.3,.3)--(6.3,.6)--(5.6,.6)--(5.6,.9)--(4.9,.9)--(4.9,1.2
      )--(4.2,1.2)--(4.2,1.5)--(3.5,1.5)--(3.5,1.8)--(2.8,1.8)--(2.8,2.1)--(2.1
      ,2.1)--(2.1,2.4)--(1.4,2.4)--(1.4,2.7)--(.7,2.7)--(.7,3)--(0,3);
  \draw[->, ultra thick,PineGreen] (0,3)--(0,7);

  \fill[TealBlue,fill=PineGreen] (7,0) circle (5pt);
  \fill[TealBlue,fill=PineGreen] (0,3) circle (5pt);
      \node at (5,4.55) (a) {\scriptsize $\widehat{A}_{10}$};
\end{tikzpicture}\quad
\begin{tikzpicture}[scale=.4]
  \path[fill=MidnightBlue!45]  (0,3)--(7,0)--(9,0)--(9,7)--(0,7)--(0,3);
  \makegrid

  \draw[->, ultra thick,MidnightBlue] (7,0)--(9,0);
  \draw[-, cap=round,ultra thick,MidnightBlue] (0,3)--(7,0);
  \draw[->, ultra thick,MidnightBlue] (0,3)--(0,7);

  \fill[TealBlue,fill=MidnightBlue] (7,0) circle (5pt);
  \fill[TealBlue,fill=MidnightBlue] (0,3) circle (5pt);
      \node at (3.5,4.55) (a) {\scriptsize $L=\lim\{\widehat{A}_p\}$};
\end{tikzpicture}\quad
\end{figure}

\end{example}

\begin{proof}
	As the origin is contained in $V$ we know that $B_p\subseteq B_p+V$ for all $p$. Further, since $B_p\subseteq L$ for all $p$ it follows that $B_p+V\subseteq L+V$ for all $p$. Thus we have that $B_p\subseteq B_p+V \subseteq L+V$ for all $p$, so parts (7) and (8) of Lemma~\ref{lem:region-limit-properties} imply
	\[
	L\subseteq \lcl \{B_p+V\} \subseteq \ucl \{B_p+V\} \subseteq L+V.
	\]
	In order to conclude the claim we show that $L+V\subseteq \lcl\{B_p+V\}$. Given $x\in L+V$ we may write $x$ as $x=\ell+v$ where $\ell\in L$ and $v\in V$. The desired claim now follows from the fact that
		\begin{eqnarray*}
			d(\ell,B_p) &=&\inf \left\{ d(\ell,y) \;\middle|\; y\in B_p \right\} = \inf \left\{ d(\ell+v,y+v) \;\middle|\; y\in B_p \right\} \\
			&\geq& \inf \left\{ d(x,z) \;\middle|\; z\in B_p+V \right\}=d(x, B_p+V).
		\end{eqnarray*}
		Since $\lim_{p\to \infty} d(\ell,B_p)=0$ it follows that $\lim_{p\to \infty} d(x, B_p+V)=0$.
\end{proof}

\begin{remark}
The other natural notion of convergence of subsets of $\RR^r$ is in the Hausdorff metric.  Since multigraded regularity is never compact this Hausdorff convergence is ill-suited for our applications.
\end{remark}

We shall prove Theorem~\ref{main:thm-ideals} in the next two sections by showing the following inclusions:
\[\ucl_{p\to\infty}\tfrac{\regR{\cI^p}}{p}\subseteq\sesh(\cI)\subseteq\lcl_{p\to\infty}\tfrac{\regR{\cI^p}}{p}\]
for the upper and lower closed limits defined in Definition~\ref{def:limits}.  The first inclusion appears in Proposition~\ref{prop:ucl-reg} and the second in Proposition~\ref{prop:lcl-reg}.

\section{Upper Closed Limit, Seshadri Regions \& Global Generation}\label{sec:ucl-sesh-global}

The goal of this section is to prove that $\ucl_{p\to\infty}\tfrac{\regR{\cI^p}}{p}$ is contained in the Seshadri region $\sesh(\cI)$. This inclusion will follow from the close connection between regularity and global generation. In particular, by \cite{MS04}*{Theorem~1.4} if $\cF$ is $D$-regular then $\cF(D)$ is globally generated. With this in mind we introduce the following set.

\begin{definition}
  Let $X$ be a smooth projective toric variety and $\cF$ an $\O_{X}$-module.
  The globally generated region of $\cF$ is
\[
  \ggen{\cF}\coloneqq \left\{ D \in \Pic(X) \;\middle|\; \cF(D) \text{ is globally generated}.\right\}
\]
  We will also need a real version:
  \begin{equation}\label{ggen}
  \ggenR{\cF}\coloneqq\ggen{\cF}+\NefR{X}.
  \end{equation}
\end{definition}

The global generation region is closed under translation by $\Nef(X)$, as in part~\eqref{it:translation} of Lemma~\ref{lem:sesh-basic-properties}, motivating the addition of $\NefR{X}$ in \eqref{ggen}.

\begin{lemma}
  Let $X$ be a smooth projective toric variety.
  If $\cF$ is an $\O_{X}$-module then
  \[\ggen{\cF}+\Nef(X)\subseteq \ggen{\cF}.\]
\end{lemma}

\begin{proof}
  The proof follows by combining two standard facts: nef line bundles are globally generated on projective toric varieties \cite{CLS11}*{Theorem 6.3.12} and products of globally generated line bundles are globally generated.  If $D\in\ggen{\cF}$ and $L\in\Nef(X)$ then $\cF(D)$ and $\O(L)$ are globally generated so $\cF(D+L)=\cF(D)\otimes\O(L)$ is as well.
\end{proof}

\subsection{Global Generation \& The Seshadri Region}

Our first lemma is a multigraded generalization of standard facts about blowups.  We use the notation $\pi^{-1}(\cI)\cdot\O_{\widetilde X}$ for the ideal sheaf in $\O_{\widetilde X}$ generated by the image of $\pi^{-1}(\cI)$, as in \cite{Hartshorne77}*{Caution 7.12.2}.

\begin{lemma}\label{lem:inverse-image-ideal-sheaves}
  Let $X$ be a smooth projective variety. Let $\cI \subseteq \O_{X}$ be an ideal sheaf and $\pi\colon\BlX \to X$ the blowup of $X$ along $\cI$ with exceptional divisor $E$. For all divisors $D\in \Pic(X)$:
  \begin{enumerate}
  \item\label{it:inverse-image} $\pi^{-1}\left(\cI^p(D)\right)\cdot\O_{\BlX}\left(\pi^*(D)\right)\cong \O_{\BlX}(\pi^*(D)-pE)$, and
  \item if $\cI^{p}(D)$ is globally generated then $\pi^{-1}\left(\cI^p(D)\right)\cdot\O_{\BlX}\left(\pi^*(D)\right)$ is globally generated.
  \end{enumerate}
\end{lemma}

\begin{proof}
	Recall that $\pi^{-1}(\mathcal I)\cdot\O_{\BlX}$ is the locally invertible sheaf $\O(-E)$ on $\BlX$ (see e.g.\ \cite{Hartshorne77}*{Proposition 7.13(a)}).  On the left hand side of \eqref{it:inverse-image} we have the image of $\pi^*\left(\mathcal I^p(D)\right)$ in $O_{\BlX}\left(\pi^*(D)\right)$.  The map
	$\mathcal I\otimes_{\O_X}\cdots\otimes_{\O_X}\mathcal I\otimes_{\O_X}\O_X(D)\to\mathcal I^p(D)$ induces a map
	\[
		\pi^*(\mathcal I)\otimes_{\O_{\BlX}}\cdots\otimes_{\O_{\BlX}}\pi^*(\mathcal I)\otimes_{\O_{\BlX}}\pi^*\left(\O_X(D)\right)\\
		\to \pi^*(\mathcal I^p(D))\to\O_{\BlX}\left(\pi^*(D)\right),
	\]
	where the image of the composition is still $\pi^{-1}\left(\cI^p(D)\right)\cdot\O_{\BlX}\left(\pi^*(D)\right)$.  The image of each $\pi^*(\mathcal I)$ in $\O_{\BlX}$ is $\O(-E)$ so we obtain a map $\O_{\BlX}(\pi^*(D)-pE)\to\pi^{-1}\left(\cI^p(D)\right)\cdot\O_{\BlX}\left(\pi^*(D)\right)$.  Locally this is an isomorphism by \cite{Vakil25}*{Exercise 15.6.E} and thus claim \eqref{it:inverse-image} follows.

	If $\cI^p(D)$ is globally generated then there is a surjection $\O_X^{\oplus n}\twoheadrightarrow\cI^p(D)$ for some $n$. By right exactness of pullback, this induces a surjection $\O_{\BlX}^{\oplus N}\twoheadrightarrow\pi^*(\cI^p(D))$.
Composing with the canonical surjection onto its image,
identified in part~(1), gives
\[
\O_{\BlX}^{\oplus N}\twoheadrightarrow
\O_{\BlX}(\pi^*(D)-pE).
\]
Thus $\O_{\BlX}(\pi^*(D)-pE)$ is globally generated.
\end{proof}

We now show that the global generation region contains the multigraded regularity when both are considered with real coefficients, a variant of \cite{MS04}*{Theorem~1.4}.

\begin{proposition}\label{prop:gg-region-regularity}
With the set-up as in Convention~\ref{notation:blowup}, if $\cF$ is a coherent sheaf on $X$ then
  \[
  \regR{\cF} \subseteq \ggenR{\cF}.
  \]
\end{proposition}
\begin{proof}
  By \cite{MS04}*{Theorem~1.4} we have $\reg(\cF)\subseteq\ggen{\cF}$, so
  \[\regR{\cF}=\reg(\cF)+\NefR{X}\subseteq\ggen{\cF}+\NefR{X}=\ggenR{\cF}\]
  according to Equations~\eqref{regR} and \eqref{ggen}.
\end{proof}

These facts imply that scaling $\ggenR{\cI^p}$ by $\frac{1}{p}$ produces elements of $\sesh(\cI)$.

\begin{proposition}\label{prop:gg-region-key}
With the set-up as in Convention~\ref{notation:blowup}, if $\cI\subseteq \O_{X}$ is an ideal sheaf then for all $p\geq 1$
  \[
  \tfrac{\ggenR{\cI^{p}}}{p} \subseteq \sesh\left(\cI\right).
  \]
\end{proposition}

\begin{proof}
  Suppose $\frac{D'}{p}\in \frac{\ggenR{\cI^{p}}}{p}$.  By \eqref{ggen} we may write $D'\in\ggenR{\cI^p}$ as $D'=D+L$ where $L\in\NefR{X}$ and $\cI^{p}(D)$ is globally generated.  Lemma~\ref{lem:inverse-image-ideal-sheaves} yields
  \[
  \pi^{-1}\left(\cI^p(D)\right)\cdot\O_{\BlX}\left(\pi^*(D)\right)\cong \O_{\BlX}\left(\pi^*(D)-pE\right)
  \]
  and $\pi^*(D)-pE$ is globally generated since  $\cI^{p}(D)$ is globally generated. As globally generated implies nef, both $\pi^*(D)-pE$ and thus also $\pi^*\left(\frac{D}{p}\right)-E$ are nef. By definition this means $\frac{D}{p}\in \sesh(\cI)$. Since $\sesh(\cI)$ is closed under translation by $\NefR{X}$ according to part~\eqref{it:translation} of Lemma~\ref{lem:sesh-basic-properties}, this implies $\frac{D'}{p}=\frac{D}{p}+\frac{L}{p}\in \sesh(\cI)$.
\end{proof}

\subsection{The Upper Closed Limit}

We now use the results from the previous subsection to establish a chain of inclusions in the limit.

\begin{proposition}\label{prop:ucl-reg}
  With the set-up as in Convention~\ref{notation:blowup},
  \[
  \ucl_{p\to\infty} \tfrac{\regR{\cI^p}}{p} \subseteq \ucl_{p\to\infty}\tfrac{\ggenR{\cI^p}}{p} \subseteq \sesh(\cI).
  \]
\end{proposition}

\begin{proof}
  We will show that for all $p\geq1$
\begin{equation}\label{eq:prop:ucl-reg}
  \tfrac{\regR{\cI^p}}{p} \subseteq \tfrac{\ggenR{\cI^p}}{p} \subseteq \sesh(\cI).
\end{equation}
  Proposition~\ref{prop:gg-region-regularity} implies that $\regR{\cI^p}\subseteq\ggenR{\cI^p}$ for all $p\geq1$. This shows the inclusion on the left side. The inclusion on the right side is Proposition~\ref{prop:gg-region-key}. Since $\sesh(\cI)$ is constant over $p$ the upper closed limit is $\overline{\sesh(\cI)}$ by part~(5) of Lemma~\ref{lem:region-limit-properties}.  By Proposition~\ref{prop:sesh-topology} the closure is unnecessary.  Thus the inclusion of limits follows from Equation~\eqref{eq:prop:ucl-reg} and part~(6) of Lemma~\ref{lem:region-limit-properties}.
\end{proof}

\section{Lower Closed Limit, Seshadri Regions \& Fujita Vanishing}\label{sec:lcl-sesh-fujita}

In this section we show that $\sesh(\cI)$ is contained in $\lcl_{p\to\infty}\tfrac{\regR{\cI^p}}{p}$. The key insight is that powerful vanishing results like Fujita's vanishing theorem apply to pullbacks of divisors in the Seshadri region. This uniform vanishing controls the multigraded regularity.

We are interested in the asymptotic behavior of $\reg\left(\cI^p\right)$, and motivated by Propositions \ref{prop:gg-region-key} and~\ref{prop:gg-region-regularity} we would like to scale this region by $\tfrac{1}{p}$.  The theorem below approximates $\reg\left(\cI^p\right)$ for high $p$ by points $s\in\sesh(\cI)$.  In the case where $s$ is rational we need to scale it by some $\ell$ to obtain an integral translate of $\Nef(X)$ depending not directly on $p$ but rather on the next multiple of $\ell$ after $p$, represented by $(q+1)\ell$ in the notation of the theorem.

\begin{theorem}\label{thm:region-regularity}
  Adopt the notation from Convention~\ref{notation:blowup}.  Fix a rational point $L\in\Int\sesh(\cI)$. Given an integer $\ell \in \ZZ_{>0}$ such that $\ell L\in N^1(X)$ (not only $\NR{X}$), there exist $m_1\in\ZZ_{>0}$ and $p_{0}\in\ZZ_{>0}$ such that
  \[
  (q+1)\ell L+n\sum_{i=1}^{t}C_{i}+\Nef(X) \subseteq \reg\left(\cI^{p}\right)
  \]
  for all $p\geq \max\{p_{0},m_1\ell\}$ where $p=q\ell+r$ for $0\leq r<\ell$.
\end{theorem}

\subsection{Fujita's Vanishing \& The Seshadri Region}

In order to prove Theorem~\ref{thm:region-regularity} we need the following two general results.

\begin{theorem}\cite{Lazarsfeld04I}*{Theorem 1.4.35, Fujita's Vanishing Theorem}\label{thm:FVT}
  Let $X$ be an irreducible projective variety, and $\cF$ a coherent sheaf on $X$. If $A$ is an ample divisor on $X$ then there exists a constant $m_0$ depending only on $\cF$ and $A$ such that for all $i>0$
  \[
  H^i(X,\cF(D))=0 \quad  \text{ for all }\quad D\in m_0A+\Nef(X).
  \]
\end{theorem}

\begin{lemma}\cite{CEL01}*{Lemma~3.3}\label{lem:blowup-cohomology-iso}
  Let $X$ be a projective variety and $\cI\subseteq \O_{X}$ an ideal sheaf. Let $\pi\colon\BlX \to X$ be the blowup of $X$ along $\cI$ with exceptional divisor $E$.
  There exists a constant $p_0\in\ZZ_{>0}$ such that for all $p\geq p_0$, all $i\geq0$, and all divisors $D$ on $X$ there is an isomorphism:
  \[
  \begin{tikzcd}
    H^i\left(X,\cI^{p}(D)\right) \rar[leftrightarrow]{\sim} &H^i\left(\BlX, \O_{\BlX}(\pi^*(D)-pE)\right).
  \end{tikzcd}
  \]
\end{lemma}

Now we will prove three technical lemmas.  The first relies on the definition of $\sesh(\cI)$ and algebra in $\Pic(\BlX)$.  The second combines this with Theorem~\ref{thm:FVT} and Lemma~\ref{lem:blowup-cohomology-iso} to describe the cohomology vanishing of high powers $\cI^p$ near points in $\sesh(\cI)$.  The third applies this vanishing to the definition of regularity.

\begin{lemma}\label{lem:nef-to-nef}
  Adopt the set-up from Convention~\ref{notation:blowup}. Fix a rational point $L\in\Int\sesh(\cI)$. Given~$L$ and integers $\ell,m_{0}\in \ZZ_{>0}$ such that $\ell L\in N^1(X)$, there exists $m_1\in\ZZ_{>0}$ such that
  \[
  D\in (q+1)\ell L+\Nef(X)\implies \pi^*(D)-pE-m_{0}\left(\pi^*(\ell L)-\ell E\right)\in \Nef(\BlX)
  \]
  for all $p>m_1\ell$ where $p=q\ell+k$ for $0\leq k<\ell$.
\end{lemma}

\begin{proof}
  Since $L\in\Int\sesh(\cI)$ the divisor $\pi^*(\ell L)-\ell E$ is ample by Proposition~\ref{prop:sesh-topology}.  First we choose~$m_1$ so that $(m_1-m_0)(\pi^*(\ell L)-\ell E)+(\pi^*(\ell L)-kE)$ is nef for all $0\leq k <\ell$, which is possible by \cite{Lazarsfeld04I}*{Theorem 1.2.6(iii)}.

  Let $D\in (q+1)\ell L+\Nef(X)$ and consider the divisor of interest on $\BlX$.  By rearranging we see that
  \begin{align*}
    \pi^*(D)&-pE-m_{0}\left(\pi^*(\ell L)-\ell E\right)\\
    &=\pi^*\left(D-m_{0}\ell L\right)-pE+m_{0}\ell E,\\
    \intertext{which using $p=q\ell+k$ we may write as}
    &= \pi^*\left(D-m_{0}\ell L\right)-(q-m_{0})\ell E - kE.
    \intertext{Adding and subtracting $(q-m_{0}+1)\pi^*(\ell L)+(m_1-m_0)(\pi^*(\ell L)-\ell E)$ gives}
    &=\begin{multlined}[t]
    \pi^*\left(D-m_{0}\ell L\right)-(q-m_0+1)\pi^*(\ell L)\\
    +(q-m_0)\pi^*(\ell L)-(q-m_{0})\ell E-(m_1-m_0)(\pi^*(\ell L)-\ell E)\\
    +\pi^*(\ell L) - kE+(m_1-m_0)(\pi^*(\ell L)-\ell E)
    \end{multlined}\\
    \intertext{which we may simplify to}
    &=\underbrace{\pi^*\left(D-(q+1)\ell L\right)}_{\text{I}}+\underbrace{(q-m_1)\left(\pi^*(\ell L)-\ell E\right)}_{\text{II}}+\underbrace{\left(\pi^*(\ell L)-kE\right)+(m_1-m_0)(\pi^*(\ell L)-\ell E)}_{\text{III}}.
  \end{align*}
  We now show that I, II, and III, are all nef divisors on $\BlX$.

  For I, since $D\in (q+1)\ell L+\Nef(X)$ we have that $D-(q+1)\ell L$ is nef on $X$. Nef divisors pull back along proper morphisms and $\pi$ is proper so $\pi^*(D-(q+1)\ell L)$ is nef on $\BlX$.

  We assumed that $p>m_1\ell$ so we know that $q\geq m_1$ and thus $q-m_1\geq 0$. We conclude that~II is nef because $L\in\sesh(\cI)$ means that $\pi^*(L)-E$ is nef by definition so positive multiples of it are as well.

  Divisor III is nef by choice of $m_1$.
\end{proof}

\begin{lemma}\label{lem:mg-ideal-vanishing}
  Adopt the set-up from Convention~\ref{notation:blowup}. Fix a rational point $L\in\Int\sesh(\cI)$. Given~$L$ and an integer $\ell \in \ZZ_{>0}$ such that $\ell L\in N^1(X)$, there exist constants $m_1$ and $p_0$ (independent of $p$ and $D$) such that
  \[
  H^{i}\left(X, \cI^{p}(D)\right) = 0 \quad \text{for all} \quad D\in (q+1)\ell L+\Nef(X),
  \]
  all $i>0$, and all $p>\max\{p_{0},m_1\ell\}$ where $p=q\ell+k$ for $0\leq k < \ell$.
\end{lemma}

\begin{proof}
  Lemma~\ref{lem:blowup-cohomology-iso} states there exists a constant $p_{0}$ depending on $\cI$ such that for all $p\geq p_{0}$ and all $i>0$ there is an isomorphism
  \begin{equation}\label{eq:lem:ideal-power-iso}
    \begin{tikzcd}
      H^i\left(X,\cI^{p}(D)\right) \rar[leftrightarrow]{\sim} &H^i\left(\BlX, \O_{\BlX}(\pi^*(D)-pE)\right).
    \end{tikzcd}
  \end{equation}
  We will show that if $p\gg0$ (in a way we will make precise) and $D\in (q+1)\ell L+\Nef(X)$ then $H^i\left(\BlX, \O_{\BlX}(\pi^*(D)-pE)\right)=0$ for all $i>0$.

  By Proposition~\ref{prop:sesh-topology}, the assumption $L\in\Int\sesh(\cI)$ implies that $\pi^*(L)-E$ is an ample divisor on $\BlX$. Therefore $\ell(\pi^*(L)-E)=\pi^*(\ell L)-\ell E$ is also ample, and further a $\ZZ$-divisor on $\BlX$ by choice of $\ell$. Therefore we may apply Fujita's Vanishing Theorem (see Theorem~\ref{thm:FVT}) for $A=\pi^*(\ell L)-\ell E$ and $\cF=\O_{\BlX}$. In particular, there exists a constant $m_{0}$ depending only on $\pi^*(\ell L)-\ell E$ and $\O_{\BlX}$ such that for all $i>0$
  \begin{equation}\label{eq:lem:FVT-ideal}
    H^{i}\left(\BlX, \O_{\BlX}(D')\right) = 0  \quad \text{for all} \quad D'\in m_{0}\left(\pi^*(\ell L)-\ell E \right)+\Nef(\BlX).
  \end{equation}
  Now by Lemma~\ref{lem:nef-to-nef} there exists $m_1\in\ZZ_{>0}$ such that for $p>m_1\ell$ and $p=q\ell+k$ where $0\leq k < \ell$
  \[
  D\in (q+1)\ell L+\Nef(X)\implies \pi^*(D)-pE-m_{0}\left(\pi^*(\ell L)-\ell E\right)\in \Nef(\BlX).
  \]
  Phrased differently, if $D\in (q+1)\ell L+\Nef(X)$ then $\pi^*(D)-pE \in m_{0}\left(\pi^*(\ell L)-\ell E\right)+\Nef(\BlX)$. This implies the desired vanishing for $p>\max\{p_0,m_1\ell\}$ by combining \eqref{eq:lem:ideal-power-iso} and \eqref{eq:lem:FVT-ideal} with $D'=\pi^*(D)-pE$.
\end{proof}

\begin{lemma}\label{lem:vanishing-to-regularity}
  Let $X$ be a smooth projective toric variety of dimension $n$ and $\cF$ a coherent sheaf on $X$. Fix $\bC=(C_{1},\ldots,C_{t})$ an ordered tuple of divisors on $X$ that are a minimal generating set for $\Nef(X)$, and fix a divisor $A$ on $X$.  If $d$ is an integer such that
  \[H^i\left(X, \cF(D)\right) = 0 \quad\text{for all}\quad D\in dA+\Nef(X)\]
  and all $i>0$ then
  \[
  dA+n\sum_{i=1}^tC_i+\Nef(X) \subseteq \reg(\cF).
  \]
\end{lemma}

\begin{proof}
  Recall that $\cF$ is $D$-regular if and only if for all $i>0$
  \[
  H^{i}\left(X, \cF(D')\right)=0 \quad \quad \text{for all} \quad \quad D' \in \bigcup_{|s|=i}\left(D-(s\cdot\bC)+\Nef(X)\right).
  \]
  Combining the fact that $H^{i}\left(X, \cF(D')\right)=0$ for all $i>n$ and the fact that
  \begin{align*}
    \bigcup_{|s|=1}\left(D-(s\cdot\bC)+\Nef(X)\right)&\subseteq \bigcup_{|s|=2}\left(D-(s\cdot\bC)+\Nef(X)\right)\subseteq \cdots \\
    &\subseteq \bigcup_{|s|=n}\left(D-(s\cdot\bC)+\Nef(X)\right)
  \end{align*}
  we see it is enough that for all $i>0$
  \[
  H^{i}\left(X, \cF(D')\right)=0 \quad \quad \text{for all} \quad \quad D' \in \bigcup_{|s|=n}\left(D-(s\cdot\bC)+\Nef(X)\right).
  \]
  The result follows by noting that if $D\in dA+n(C_{1}+\cdots+C_{t})+\Nef(X)$ and $|s|=n$ then $D-(s\cdot\bC)$ is in $dA+\Nef(X)$.
\end{proof}

\begin{proof}[Proof of Theorem~\ref{thm:region-regularity}]
  The desired claim follows by combining Lemmas \ref{lem:mg-ideal-vanishing} and~\ref{lem:vanishing-to-regularity}.
\end{proof}

\subsection{The Lower Closed Limit}

We know from Proposition~\ref{prop:gg-region-regularity} that the regularity of an ideal sheaf $\cI$ is contained in its global generation region, and we saw in Proposition~\ref{prop:gg-region-key} that the contraction by $\tfrac{1}{p}$ of the global generation region of $\cI^p$ is contained in the Seshadri region of $\cI$.  Therefore we may compare the contraction by $\tfrac{1}{p}$ of the regularity of $\cI^p$ to the Seshadri region of $\cI$.

\begin{proposition}\label{prop:lcl-reg}
  With notation as in Convention~\ref{notation:blowup},
  \[
  \sesh(\cI)\subseteq \lcl_{p\to\infty} \tfrac{\regR{\cI^p}}{p}.
  \]
\end{proposition}

\begin{proof}
  For convenience let $B_p=\tfrac{\regR{\cI^p}}{p}$ for each $p>0$.  Since $\sesh(\mathcal I)=\overline{\Int\sesh(\mathcal I)}$ by Proposition~\ref{prop:sesh-topology} and $\lcl B_p$ is closed By Lemma~\ref{lem:region-limit-properties}, it suffices to prove the inclusion
  \[\Int\sesh(\cI)\subseteq \lcl_{p\to\infty} B_p.\]

  We wish to show that if $x\in\Int\sesh(\cI)$ then $\limsup \left\{d(x,B_p)\right\}=0$. Applying the characterization of lower closed limits from Lemma~\ref{lem:unpacking-region-limit-defs} this is equivalent to showing that for all $\epsilon>0$ there exists $N>0$ such that $d(x,B_p)\leq \epsilon$ for all $p\geq N$. However, given a different point $s\in\Int\sesh(\cI)$ the triangle inequality implies that
  \[
  d\left(x,B_p\right)\leq d\left(x,s\right)+d\left(s,B_p\right)
  \]
  for all $p$. Therefore it is enough to show: for all $\epsilon>0$ there exists a rational point $s\in\Int\sesh(\cI)$ such that $d(x,s)\leq \frac{\epsilon}{2}$ and there exists $N>0$ such that $d(s,B_p)\leq \frac{\epsilon}{2}$ for all $p\geq N$. We will actually show the stronger claim that the second condition is true for all rational points $s\in\Int\sesh(\cI)$ such that $d(x,s)\leq \frac{\epsilon}{2}$ (although the constant $N$ depends on $s$).

  Fix $\epsilon>0$. Since $\sesh(\cI)$ is the closure of its interior by  Proposition~\ref{prop:sesh-topology}, we can choose a rational point $s\in\Int\sesh(\cI)$ such that $d(x,s)\leq \frac{\epsilon}{2}$.

  Let $\ell\in \ZZ_{>0}$ such that $\ell s\in N^1(X)$. By Theorem~\ref{thm:region-regularity} we know if $p>\max\{p_{0},m_1\ell\}$ then $(q+1)\ell s+n\sum_{i=1}^t C_{i}$ is in the regularity of $\cI^{p}$ where  $p=q\ell+k$ with $0\leq k <\ell$. In particular, consider the point
  \begin{align*}
    s_p\coloneqq\frac{(q+1)\ell s+n\sum_{i=1}^{t}C_{i}}{p} \in \frac{\regR{\cI^p}}{p}=B_p
  \end{align*}
  where
  \[
  s-s_p = \frac{(q\ell+k)s}{p}-\frac{(q+1)\ell s+n\sum_{i=1}^{t}C_{i}}{p}=\frac{k-\ell}{p}s-\frac{n}{p}\sum_{i=1}^t C_{i}
  \]
  Taking $p\gg0$ we can ensure that $\left| s-s_p \right |\leq \frac{\epsilon}{2}$, where the norm is taken in $\RR^r$.  This means that as needed
  \[
  d\left(s,B_p\right)\leq d\left(s,s_p\right) \leq \frac{\epsilon}{2}.
  \]
\end{proof}

\subsection{Proofs of the main theorems}\label{sub-sec:proof-of-main}

\begin{proof}[Proof of Theorem \ref{main:thm-ideals}]
The inclusions in Proposition~\ref{prop:ucl-reg} and Proposition~\ref{prop:lcl-reg} show that
	\[\ucl_{p\to\infty}\tfrac{\regR{\cI^p}}{p}\subseteq\sesh(\cI)\subseteq\lcl_{p\to\infty}\tfrac{\regR{\cI^p}}{p}.\]
The statement now follows from Lemma~\ref{lem:region-limit-properties} which shows that in general the lower closed limit is contained in the upper closed limit.
\end{proof}

\begin{proof}[Proof of Corollary~\ref{main:cor-convex}]
This follows from Theorem~\ref{main:thm-ideals} since $\sesh(\cI)$ is closed and convex (see Proposition~\ref{prop:sesh-topology}) and defined by finitely many rational hyperplanes when $\BlX$ is a Mori dream space (see Corollary~\ref{cor:sesh-poly}).
\end{proof}

\section{A parallel theorem for vector bundles}\label{sec:parallel-bundles}

In this section we adapt the techniques developed in the previous sections to describe the asymptotic regularity of symmetric powers of a vector bundle $\mathcal E$. The arguments closely parallel the ideal sheaf case, relying on an analogy between blowups and projective bundles.

Following the conventions in \cite{Hartshorne77}*{Chapter II.7} we let $\Sym^{\bullet}(\cE) \coloneqq \bigoplus_{p\geq0}\Sym^p(\cE)$ be the symmetric algebra of $\cE$, and define the projective bundle $\PP(\cE)$ to be $\Proj(\Sym^{\bullet}(\cE))$. Let $\pi\colon\PP(\cE)\to X$ be the natural projection morphism and $\O_{\PP(\cE)}(1)$ the invertible sheaf that arises from this construction.  For this section, let $F$ be a Cartier divisor corresponding to $\O_{\PP(\cE)}(1)$.

Like the cohomology on a blowup, the cohomological vanishing of $\Sym^{p}(\cE)$ can be translated to vanishing questions for $\O_{\PP(\cE)}(p)$, which we can again deal with uniformly. In fact, in many ways the vector bundle case is more natural.  We define the Seshadri region of $\cE$ in this situation as follows.

\begin{definition}\label{def:seshBundles}
  Let $X$ be a projective variety and $\cE$ a vector bundle on $X$ with $\O_{\PP(\cE)}(1)=\O_{\PP(\cE)}(F)$.  The \emph{Seshadri region of $\cE$ restricted to} a subspace $V \subset\NR{X}$ is
  \[
  \sesh_V\left(\cE\right) \coloneqq \left\{L\in V \; \middle| \; \pi^*(L)+F \text{ is nef on $\PP(\cE)$} \right\} \subset\NR{X}.
  \]
  The \emph{Seshadri region} of $\cE$ is the Seshadri region of $\cE$ restricted to the $\RR$-span of $\NefR{X}$ and we denote this region by $\sesh(\cE)\coloneqq \sesh_{\NefR{X}}(\cE)$.
\end{definition}

\begin{remark}
  At first glance the presence of $+F$ instead of $-F$ may be surprising since Definition~\ref{def:sesh} has $-E$. This confusion arises from our somewhat unfortunate combination of standard conventions for divisors and line bundles.  As in the proof of Lemma~\ref{lem:inverse-image-ideal-sheaves}, the inverse image ideal sheaf $\pi^{-1}\cI \cdot \O_{\BlX}$ in the blowup case is the same as $\O_{\BlX}(-E)$.
\end{remark}

\begin{remark}
	As a warning, the Seshadri region of a vector bundle $\cE$ defined above is different from the local and global Seshadri constants of vector bundles in \cites{BSS93,BSS96,hacon00,HMP10}. Those constants are defined using blowups at points: the local Seshadri constant at $x\in X$ is the supremum over $\lambda\in \QQ_{>0}$ such that $\O_{\PP(\pi^*(\cE))}(1)-\lambda\widetilde{E}$ is nef on $\PP(\pi^*(\cE))$ where $\pi\colon\BlX\to X$ is the blowup at $x$ with exceptional divisor $E$, and $\widetilde{E}$ is the pull-back of $E$ to the projectivization of $\pi^*(\cE)$. The global Seshadri constant is the infimum over all $x\in X$. In contrast our Seshadri regions do not involve blowups of $X$.
\end{remark}

One can show that almost everything proven in Section~\ref{sec:seshadri-region} concerning Seshadri regions of ideal sheaves is true, verbatim, for Seshadri regions of vector bundles. We record one basic lemma which describes nef divisors on $\PP(\cE)$ in terms of nef vector bundles.  A vector bundle is nef when $\O(1)$ is nef on its projectivization.

\begin{lemma}\label{lem:vb-basic}
	Let $X$ be a projective variety, $\cE$ a vector bundle on $X$, and $\pi\colon\PP(\cE)\to X$ the natural projection map. A divisor $L \in N^1(X)$ is in $\sesh(\cE)$ if and only if $\cE(L)$ is nef.
\end{lemma}

\begin{proof}
  By \cite{Hartshorne77}*{Lemma~7.9} there is an isomorphism between $\PP(\cE(L))$ and $\PP(\cE(L))$ sending $\O_{\PP(\cE(L))}(1)$ to $\O_{\PP(\cE)}(\pi^*(L)+F)$, so one is nef if and only if the other is.
\end{proof}

\begin{example}
	When restricted to a line the tangent bundle $T_{\PP^{n}}$ of $\PP^n$ is isomorphic to $\O_{\PP^1}(2)\oplus \O_{\PP^1}(1)^{n-1}$ for $n\geq 2$.  If $T_{\PP^n}$ is nef then this restriction is as well, so $\sesh(T_{\PP^{n}}) \subset \{aH \; | \; a \geq-1\}$ where $H$ is the hyperplane class on $\PP^n$.  Since $T_{\PP^{n}}(-1)$ is a quotient of a globally generated bundle via the Euler sequence it is itself globally generated, hence nef.  Thus the Seshadri region of $T_{\PP^{n}}$ is  $\{ aH \; | \; a \geq-1\}$. A similar argument shows that $\sesh(T_{\PP^1\times\PP^1}) = \{ aH_1 + bH_2 \; | \; a,b\geq0\}$ where $H_1$ and $H_2$ are the fibers of the projections to $\PP^1$.
\end{example}

Analogously to our results for powers of ideal sheaves, the asymptotic regularity of symmetric powers of a vector bundle is equal to its Seshadri region.

\begin{theoremalpha}\label{main:thm-bundles}
  Let $X$ be a smooth projective toric variety. If  $\cE$ is a locally-free sheaf on $X$ then
  \[
  \lim_{p\to\infty} \tfrac{\regR{\Sym^p(\cE)}}{p} =  \sesh\left( \cE \right).
  \]
\end{theoremalpha}

The proof of Theorem~\ref{main:thm-bundles} largely mimics that of Theorem~\ref{main:thm-ideals}, with the exception of the following somewhat standard lemma.

\begin{lemma}\label{lem:projectiviztion-cohomology-iso}
  Let $X$ be a projective variety and $\cE$ a vector bundle on $X$ with $\rank \cE \geq2$. For all $p>0$, all $i\geq0$, and all divisors $D$ on $X$ there is an isomorphism:
  \[
  \begin{tikzcd}
    H^i\left(X,\Sym^p(\cE)\otimes\O_X(D)\right) \rar[leftrightarrow]{\sim} &H^i\left(\PP(\cE), \pi^*\O_{X}(D)\otimes \O_{\PP(\cE)}(p)\right)
  \end{tikzcd}
  \]
\end{lemma}

\begin{proof}
This amounts to stringing together a series of exercises from Section III.8 of \cite{Hartshorne77}. By \cite{Hartshorne77}*{Exercise III.8.4} we know that $R^i\pi_*\O_{\PP(\cE)}(p)=0$ for all $i>0$ and all $p>0$. Combining this with the projection formula (see \cite{Hartshorne77}*{Exercise III.8.3}) gives
\[
R^i\pi_*\left(\O_{\PP(\cE)}(p) \otimes \pi^*\O_{X}(D)\right) \cong R^i\pi_*\left(\O_{\PP(\cE)}(p)\right) \otimes \O_{X}(D) = 0
\]
for all $i>0$, all divisors $D$ on $X$, and all $p>0$. This vanishing allows us to invoke the degenerate version of the Leray spectral sequence (see \cite{Hartshorne77}*{Exercise III.8.1}) which implies that for all divisors $D$ and all integers $p>0$
 \[
  H^i\left(X, \pi_*\left(\O_{\PP(\cE)}(p) \otimes \pi^*\O_X(D)\right)\right) \cong H^i\left(\PP(\cE), \O_{\PP(\cE)}(p)\otimes \pi^*\O_{X}(D) \right).
  \]
The claimed result then follows from the fact that $\pi_*\O_{\PP(\cE)}(p)\cong \Sym^{p}(\cE)$ for all $p>0$  \cite{Hartshorne77}*{Proposition II.7.11.(a)}.
\end{proof}

In addition to the above lemma the other input needed for the proof of Theorem~\ref{main:thm-bundles} is an analogue of the vanishing results given in Section~\ref{sec:lcl-sesh-fujita}.

\begin{proposition}
	Let $X$ be a smooth projective toric variety of dimension $n$ and $\cE$ a vector bundle on $X$ with $\rank \cE \geq2$ and $\O_{\PP(\cE)}(1)=\O_{\PP(\cE)}(F)$. Fix a rational point $L \in \Int\sesh(\cE)$ and $\ell \in \ZZ_{>0}$ such that $\ell L \in  N^1(X)$.
	\begin{enumerate}[(i)]
    \item There exists a constant $m_{0}$ depending on $\cE$ such that for all $p>m_0\ell$ and all $i>0$
    \[
    H^{i}\left(X, \Sym^{p}(\cE)\otimes\O_X (D)\right) = 0 \quad \text{for all} \quad D\in  \left(\left \lfloor \tfrac{p}{\ell}\right \rfloor+1\right)\ell  L+\Nef(X).
    \]
    \item Given $m_{0}\in \ZZ_{>0}$ there exists $m_1\in\ZZ_{>0}$ such that for all $p>m_1\ell$
    \[
    D\in \left(\left \lfloor \tfrac{p}{\ell}\right \rfloor+1\right)\ell L+\Nef(X)\implies \pi^*(D)+pF-m_{0}\left(\pi^*(\ell L)+\ell F\right)\in \Nef(\PP(\cE)).
    \]
    \item Let $C_{1},\ldots,C_{t}$ generate $\Nef(X)$ and let $m_1$ be the constant from (ii). For all $p\geq m_1\ell$
    \[
    \left(\left \lfloor \tfrac{p}{\ell}\right \rfloor+1\right)\ell L+n\sum_{i=1}^{t}C_{i}+\Nef(X) \subseteq \reg\left(\Sym^p(\cE)\right).
    \]
	\end{enumerate}
\end{proposition}

The proofs are analogous to those in Section~\ref{sec:lcl-sesh-fujita}, so we suppress them here.  Specifically part~(i) is Lemma~\ref{lem:mg-ideal-vanishing}, part~(ii) is Lemma~\ref{lem:nef-to-nef}, and part~(iii) is Theorem~\ref{thm:region-regularity}.

\section{Example: blowing up the del Pezzo surface of degree 7}\label{sec:examples}

Consider the toric variety $\PP^1\times\PP^1$ with coordinates $x_0,x_2$ on one factor and $x_1,x_3$ on the other.  The blowup of $\PP^1\times\PP^1$ at the torus-invariant point $V(x_1,x_2)$ is a smooth projective toric variety which appeared in Example~\ref{ex:blowup-P1xP1}. In fact, it is the Fano variety $F_{2,3}$. The fan associated to $F_{2,3}$ is the complete two-dimensional fan with the five rays $\rho_0=(1,0)$, $\rho_1=(0,1)$, $\rho_2=(-1,1)$, $\rho_3=(-1, 0)$, and $\rho_4=(0, -1)$, shown in Figure \ref{fig:fans} below.

Let $S=\CC[x_0,x_1,x_2,x_3,x_4]$ be the homogeneous coordinate ring of $F_{2,3}$.  (Note that these coordinates do not agree with the coordinates on $\PP^1\times\PP^1$ above, but rather are chosen to agree with the indexing in \emph{Macaulay2}.)  We compute the Seshadri region of the ideal sheaf corresponding to the ideal $I=\langle x_0,x_4\rangle$. Blowing up $F_{2,3}$ at the torus invariant point where $I$ vanishes, one obtains another smooth toric Fano variety $F_{2,4}$. The fan associated to $F_{2,4}$ contains all the rays of the previous fan as well as the additional ray $\rho_5=(1,-1)$, as shown in Figure~\ref{fig:fans}.

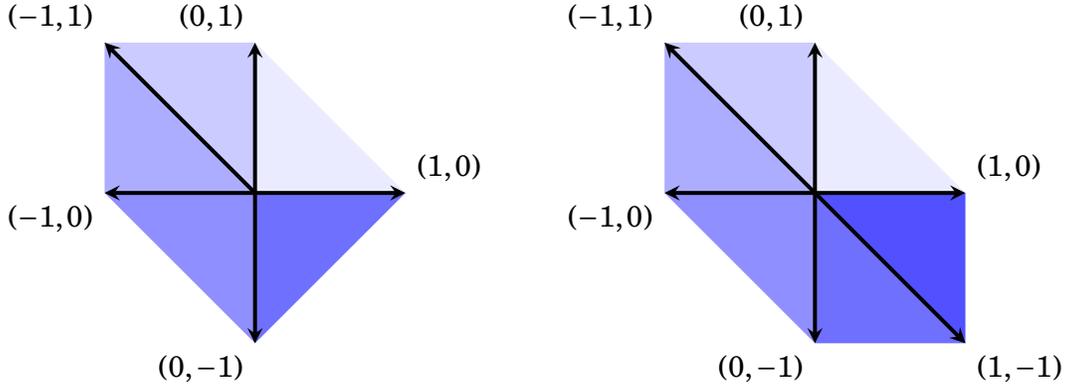
\begin{figure}[h!]

  \begin{tikzpicture}[scale=2, >=stealth]
    \coordinate (V0) at (1,0);
    \coordinate (V1) at (0,1);
    \coordinate (V2) at (-1,1);
    \coordinate (V3) at (-1,0);
    \coordinate (V4) at (0,-1);

    \fill[blue!10, opacity=0.8] (0,0) -- (V0) -- (V1) -- cycle;
    \fill[blue!25, opacity=0.8] (0,0) -- (V1) -- (V2) -- cycle;
    \fill[blue!40, opacity=0.8] (0,0) -- (V2) -- (V3) -- cycle;
    \fill[blue!55, opacity=0.8] (0,0) -- (V3) -- (V4) -- cycle;
    \fill[blue!70, opacity=0.8] (0,0) -- (V4) -- (V0) -- cycle;

    \draw[->, line width=1.5pt,black,-stealth](0,0) -- (V0) node[anchor=south west] {$(1,0)$};
    \draw[->, line width=1.5pt,black,-stealth] (0,0) -- (V1) node[anchor=south east] {$(0,1)$};
    \draw[->, line width=1.5pt,black,-stealth] (0,0) -- (V2) node[anchor=south east] {$(-1,1)$};
    \draw[->, line width=1.5pt,black,-stealth] (0,0) -- (V3) node[anchor=north east] {$(-1,0)$};
    \draw[->, line width=1.5pt,black,-stealth] (0,0) -- (V4) node[anchor=north east] {$(0,-1)$};

  \end{tikzpicture}
  \qquad
  \begin{tikzpicture}[scale=2, >=stealth]

    \coordinate (V0) at (1,0);
    \coordinate (V1) at (0,1);
    \coordinate (V2) at (-1,1);
    \coordinate (V3) at (-1,0);
    \coordinate (V4) at (0,-1);
    \coordinate (V5) at (1,-1); 

    \fill[blue!10, opacity=0.8] (0,0) -- (V0) -- (V1) -- cycle;
    \fill[blue!25, opacity=0.8] (0,0) -- (V1) -- (V2) -- cycle;
    \fill[blue!40, opacity=0.8] (0,0) -- (V2) -- (V3) -- cycle;
    \fill[blue!55, opacity=0.8] (0,0) -- (V3) -- (V4) -- cycle;
    \fill[blue!70, opacity=0.8] (0,0) -- (V4) -- (V5) -- cycle;
    \fill[blue!85, opacity=0.8] (0,0) -- (V5) -- (V0) -- cycle;

    \draw[->, line width=1.5pt,black,-stealth]   (0,0) -- (V0) node[anchor=south west] {$(1,0)$};
    \draw[->, line width=1.5pt,black,-stealth]   (0,0) -- (V1) node[anchor=south east] {$(0,1)$};
    \draw[->, line width=1.5pt,black,-stealth] (0,0) -- (V2) node[anchor=south east] {$(-1,1)$};
    \draw[->, line width=1.5pt,black,-stealth] (0,0) -- (V3) node[anchor=north east] {$(-1,0)$};
    \draw[->, line width=1.5pt,black,-stealth] (0,0) -- (V4) node[anchor=north east] {$(0,-1)$};
    \draw[->, line width=1.5pt,black,-stealth] (0,0) -- (V5) node[anchor=north west] {$(1,-1)$};

  \end{tikzpicture}
  \caption{The fans of the Fano surfaces $F_{2,3}$ (left) and $F_{2,4}$ (right)}.
  \label{fig:fans}
\end{figure}

The Picard group of $F_{2,3}$ is isomorphic to the sub

lattice of $\ZZ^3$  generated by the columns of the following matrix, with the multidegree of each variable in $S$ given by the corresponding column:
\begin{equation}\label{eq:degree-matrix}
  \begin{bmatrix}
    1 & -1 & 1 & 0 & 0 \\
    0 & 1 & -1 & 1 & 0 \\
    0 & 0 & 1 & -1 & 1
  \end{bmatrix}.
\end{equation}
In these coordinates $\NefR{F_{2,3}}$ is identified with the positive orthant in $\RR^3$.

For the rest of the section we discuss the following assertion in the context of our results.

\begin{figure}[!ht]
  \centering
  \resizebox{0.5\textwidth}{!}{%
    \begin{tikzpicture}
      \tikzstyle{every node}=[font=\scriptsize]
       \path[fill=PineGreen!30] (10.5,15.25)--(7.50,18.50)--(7.5,15.25)--cycle;
        \path[fill=PineGreen!30] (5.00, 12.50)--(10.5,15.25)--(7.5,15.25)--cycle;
          \path[fill=PineGreen!45] (5.5,15.75) -- (7.50, 15.25)--(5.00, 12.50)-- (3.60, 14.00)--cycle;
           \path[fill=PineGreen!45] (5.5,15.75) -- (7.50, 15.25)--(7.50, 18.50)--(5.50, 18.00)--cycle;
            \path[fill=PineGreen!30] (3.60, 14.00)--(5.5,15.75) -- (5.50, 18.00)--cycle;
      \draw [->, dashed] (6.25,15.25) -- (6.25,18.50) node[above]{z};
      \draw [->, dashed] (6.25,15.25) -- (3.75,12.75) node[below]{x};
      \draw [->, dashed] (6.25,15.25) -- (10.5,15.25) node[above]{y};
      \node [] at (8,15) {(0,1,0)};
      \node [] at (4.9,16) {(1,0,1)};
      \draw [->] (5.5,15.75) -- (5.50, 18.00); 
      \draw [->] (7.5,15.25) -- (7.50, 18.50); 
      \draw [->, line width=0.5pt] (7.5,15.25) -- (10.5, 15.25); 
      \draw [->] (7.5,15.25) -- (5.00, 12.50); 
      \draw [->] (5.5,15.75) -- (3.60, 14.00); 
      \draw []   (5.5,15.75) -- (7.50, 15.25); 
    \end{tikzpicture}
  }\caption{Seshadri region of $F_{2,3}$ at $\cI$}\label{fig:Seshadri}
\end{figure}
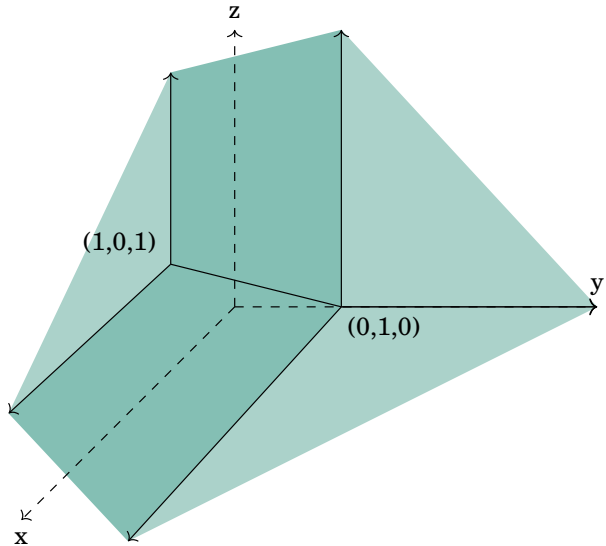

\begin{proposition}
The Seshadri region of the ideal sheaf $\cI=\widetilde{\langle x_0,x_4\rangle}$ on the Fano variety $F_{2,3}$ is the convex body \[
S(\mathcal I)
=\{(a,b,c)\in\mathbf{R}_{\ge0}^3:
       a+b\ge1,\ b+c\ge1\}
=\operatorname{conv}\{(1,0,1),(0,1,0)\}
  +\mathbf{R}_{\ge0}^3.
\]
 shown in Figure \ref{fig:Seshadri} and extending in the direction of the positive orthant.
\end{proposition}
\begin{proof}
Leveraging the fact that the blowup of $F_{2,3}$ at $\cI$ is a toric variety, one can compute the Seshadri region using Definition \ref{def:sesh} by performing linear algebra in $\PicR{F_{2,4}}$. 

Let $D_i$ denote the torus-invariant curve on $F_{2,3}$ corresponding to the ray $\rho_i$, and let $\widetilde D_i$ denote its strict transform under the blowup $\pi:F_{2,4}\to F_{2,3}$ at $D_0\cap D_4$. In the coordinates of (8.1), a class $L=(a,b,c)$ is represented by
\[
(a+b)D_0+bD_1+cD_4.
\]
Since the rays corresponding to $D_0,D_1,D_2,D_3,D_4$ occur in cyclic order in the fan, neighboring curves intersect once and non-neighboring curves are disjoint. For self-intersections, use the
primitive ray generators $v_i$: if
\[
v_{i-1}+v_{i+1}=b_i v_i,
\]
(with indices taken modulo $5$) then $D_i^2=-b_i$ \cite{CLS11}*{Theorem~10.4.4}.

Hence $(D_0^2,D_1^2,D_2^2,D_3^2,D_4^2)=(0,-1,-1,-1,0)$, so
\[
L\cdot D_0=b+c,\quad L\cdot D_1=a+b-b=a,\quad L\cdot D_2=b,\quad L\cdot D_3=c,\quad L\cdot D_4=a+b.
\]
Since the exceptional curve $E$ has self-intersection $-1$ and meets only $\widetilde D_0$ and $\widetilde D_4$, each transversely once, the intersections of $\pi^*L-E$ with $\widetilde D_0,\ldots,\widetilde D_4$ are
\[b+c-1,\quad a,\quad b,\quad c,\quad a+b-1,\]
while $(\pi^*L-E)\cdot E=1$. A divisor on a complete toric surface is nef if and only if its intersection with every invariant curve is nonnegative. Therefore the Seshadri region is \[
S(\mathcal I)
=\{(a,b,c)\in\mathbf{R}_{\ge0}^3:
       a+b\ge1,\ b+c\ge1\}
=\operatorname{conv}\{(1,0,1),(0,1,0)\}
  +\mathbf{R}_{\ge0}^3.
\]
as depicted in Figure 4.
\end{proof}

 The proof of Proposition~\ref{prop:sesh-toric} contains a similar computation, except there we did not compute explicit equations of facets in $\PicR{X}$.

Next we illustrate our Theorem \ref{main:thm-ideals} by showing how the asymptotic regularity region relates to Figure \ref{fig:Seshadri}. While in this case computing the Seshadri region was made accessible by readily available information on its nef cone, the method illustrated below may prove more tractable in situations where such information is absent.  This computation also illustrates the inherent difficulty of determining multigraded regularity regions. Indeed, very few examples of such computations are available in the literature.

\begin{proposition}
Let $S=k[x_0, x_1, x_2, x_3,x_4]$ be the Cox ring of $F_{2,3}$ and let $\cI=\widetilde{\langle x_0,x_4\rangle}$. The convex body bounded by the facets shown in Figure \ref{fig:Seshadri} and extending in the direction of the positive orthant is contained in $\lim_{p\to\infty}\tfrac{\regR{\cI^p}}{p}$.
\end{proposition}

\begin{proof}
It suffices to verify that the points $(1,0,1)$ and $(0,1,0)$ are in $\lim_{p\to\infty} \tfrac{\regR{\cI^p}}{p}$. Since by Corollary \ref{main:cor-convex} the limit regularity region is convex and closed under adding the Nef cone, which in this case is the positive orthant,  we obtain  that the region pictured in Figure \ref{fig:Seshadri} is contained in $\lim_{p\to\infty}\tfrac{\regR{\cI^p}}{p}$.

Let $B$ be the toric irrelevant ideal of $F_{2,3}$.  For $I=\langle x_0, x_4\rangle$, the powers of the ideal sheaf $\cI=\widetilde I$ are the saturated powers $\cI^m=\widetilde{(I^m):B^\infty}$. Since $I$ is generated by a regular sequence and is $B$-saturated, its powers are also $B$-saturated \cite{ZariskiSamuel}*{Lemma 5, Appendix 6}, so $\cI^m=\widetilde{(I^m)}$ and $\reg(I^p)=\reg(\cI^p)$. Since $I$ is a complete intersection, the minimal free resolution of $I^m$ is an Eagon--Northcott complex:
\begin{equation}\label{eq:EN}
  0\to \bigoplus_{a+b=m-1} S(-a-1,0,-b-1)\to  \bigoplus_{a+b=m} S(-a,0,-b)\to I^m\to 0.
\end{equation}

{\em Claim A:} For all integers $m\geq 1$, the ideal $I^m$ is regular at $(m,0,m)\in\Pic(X)$.

The argument is analogous to that in Example~\ref{ex:blowup-P1xP1}, using \cite{MS04}*{Corollary~7.3} and $
\mathbb R_{\geq0}^3\subseteq\regR{S}$:

\begin{align*}
  \regR{I^m}&\supseteq \bigcup_{\phi\colon[3]\to[3]}\left(\left[\bigcap_{a+b=m}(a,0,b)+\RR_{\geq0}^3\right]\cap\left[\bigcap_{a+b=m-1}(a+1,0,b+1)-e_{\phi(1)}+\RR_{\geq0}^3\right]\right)\\
  &=(m,0,m)+\RR_{\geq0}^3.
\end{align*}

{\em Claim B:} For all integers $m\geq 1$, the ideal $I^m$ is regular at $(0,m,0)\in\Pic(X)$.

By \eqref{eq:EN}, the resolution of $I^m(0,m,0)$ is
\[
0\to \bigoplus_{a+b=m-1} S(-a-1,m,-b-1)\to  \bigoplus_{a+b=m} S(-a,m,-b)\to I^m(0,m,0) \to 0,
\]
leading to the long exact sequence on local cohomology supported at the toric irrelevant ideal $B$ of $F_{2,3}$:
\[
\small{
  \begin{tikzcd}[column sep=small]
    0 \arrow[r] & \underbrace{\bigoplus_{a+b=m-1} H_B^0(S)(-a\!-\!1,m,-b\!-\!1)}_{\mathclap{0}} \arrow[r] &  \underbrace{\bigoplus_{a+b=m} H_B^0(S)(-a,m,-b)}_{\mathclap{0}} \arrow[r] & H_B^0(I^m)(0,m,0) \\
    {} \arrow[r] & \underbrace{\bigoplus_{a+b=m-1} H_B^1(S)(-a-1,m,-b-1)}_{\mathclap{0}} \arrow[r] & \underbrace{\bigoplus_{a+b=m} H_B^1(S)(-a,m,-b)}_{\mathclap{0}} \arrow[r] & H_B^1(I^m)(0,m,0) \\
    {} \arrow[r] & \bigoplus_{a+b=m-1} H_B^2(S)(-a-1,m,-b-1) \arrow[r, "\varphi"] & \bigoplus_{a+b=m} H_B^2(S)(-a,m,-b) \arrow[r] & H_B^2(I^m)(0,m,0) \\
    {} \arrow[r] & \bigoplus_{a+b=m-1} H_B^3(S)(-a-1,m,-b-1) \arrow[r, "\psi"]  & \bigoplus_{a+b=m} H_B^3(S)(-a,m,-b) \arrow[r] & H_B^3(I^m)(0,m,0).
  \end{tikzcd}
}
\]

Higher terms vanish because $\dim X=2$, and $H_B^0(S)=H_B^1(S)=0$ because $S$ is saturated.

To establish {\em Claim B} it suffices to show:
\begin{enumerate}
  \item\label{it:psi-source}  $\psi$ has source the zero module in degrees $(-i,-j,-k)$ with $i+j+k=1$,
  \item\label{it:psi-target} $\psi$ has target the zero module in degrees $(-i,-j,-k)$ with $i+j+k=2$,
  \item\label{it:phi-iso} $\varphi$ is an isomorphism in degrees $(-1,0,0),(0,-1,0),(0,0,-1)$ and
  \item\label{it:phi-inj} $\varphi$ is an injection in degree $(0,0,0)$.
\end{enumerate}

According to \cite{EMS00}*{Theorem 1.1} the presence of local cohomology in a given fine degree $\alpha\in\ZZ^5$ depends only on the set $\Neg\alpha=\{i \;|\;\alpha_i<0\}$.
By \cite{MS04}*{Proposition 3.2} this is determined by the reduced homology of induced subcomplexes of a simplicial complex $\Delta$ whose vertices correspond to rays in the fan.  For $F_{2,3}$ the complex $\Delta$ is a pentagon:

\begin{figure}[!ht]
  \centering
  \resizebox{0.25\textwidth}{!}{%
    \begin{tikzpicture}
      \foreach \i in {0,...,4} {
        \coordinate (P\i) at ({90 + \i*72}:3);
      }

      \draw[thick] (P0) -- (P1) -- (P2) -- (P3) -- (P4) -- cycle;

      \foreach \i in {0,...,4} {
        \node at (P\i) [circle,fill=black,inner sep=1pt,label={90 + \i*72:$x_{\i}$}] {};
      }
    \end{tikzpicture}
  }
  \caption{The simplicial complex $\Delta$ used to compute the cohomology of $F_{2,3}$}
\end{figure}
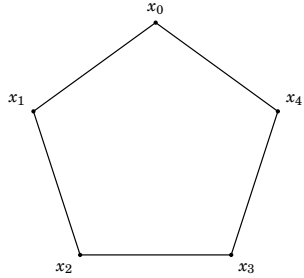

The top local cohomology module $H_B^3(S)$ has elements in degrees $(-1,-1,-1)-\Eff(X)$ where $(-1,-1,-1)$ is the degree of the canonical divisor.  This corresponds to the empty induced subcomplex of $\Delta$.

If $i+j+k=1$ then each twist appearing in degree $(-i,-j,-k)$ of the source of $\psi$ has coordinates with sum $-a-1-i+m-j-b-1-k=(m-1-a-b)-(i+j+k)-1=0-1-1=-2$.  However the sums of coordinates of degrees in $(-1,-1,-1)-\Eff(X)$ are at most $-3$, so we conclude \eqref{it:psi-source} from Claim~B.  Similarly, if $i+j+k=2$ then each twist appearing in degree $(-i,-j,-k)$ of the target of $\psi$ has coordinates with sum $-a-i+m-j-b-k=(m-a-b)-(i+j+k)=0-2=-2$, so we conclude part \eqref{it:psi-target} from Claim~B.

By \cite{MS04}*{Proposition 3.2} nonzero elements in $H_B^2(S)$ correspond to induced subcomplexes $\Delta'$ of $\Delta$ with $\widetilde{H}^0(\Delta')\neq 0$, meaning that $\Delta'$ is not connected.  The possible choices for $\Delta'$ are listed in the leftmost column of the table in Figure~\ref{fig:cohomology-types}.

To determine the source of $\varphi$ in degrees $(-1,0,0),(0,-1,0),(0,0,-1)$ we must identify Laurent monomials $x^\alpha$ with $\Neg\alpha$ equal to some $\Delta'$ and $\deg x^\alpha$ equal to $(-a-2,m,-b-1)$, $(-a-1,m-1,-b-1)$, or $(-a-1,m,-b-2)$ for $a+b=m-1$.  Similarly $x^\alpha$ in the target of $\varphi$ must have degree equal to $(-a-1,m,-b)$, $(-a,m-1,-b)$, or $(-a,m,-b-1)$ for $a+b=m$.

Recall that $\deg x^\alpha$ is the product of the matrix \eqref{eq:degree-matrix} with $\alpha$.  For most choices of $\Delta'$ the monomials $x^\alpha$ with $\Neg\alpha=\Delta'$ have degrees restricted to values which do not appear in the source or target of $\varphi$.  As an example, if $\Neg\alpha=\{1,3\}$ then multiplying the matrix \eqref{eq:degree-matrix} by $\alpha$ can only give positive values in the first and third coordinates of $\deg x^\alpha$.  Since $a,b\geq 0$ in the degrees listed above, none of them have this property.  Similar reasons for the exclusion of other possible $\Delta'$ appear in the rightmost column of the table in Figure~\ref{fig:cohomology-types}.

\begin{figure}
  \renewcommand{\arraystretch}{1.5}
  \begin{tabular}{l|l}
    $\Neg\alpha=\Delta'$ & presence of $x^\alpha$ in source and target of $\varphi$\\\hline
    $\{0,2\}$ & present in both source and target\\
    $\{1,3\}$ & absent because the first and last coordinates of $\deg x^\alpha$ must be positive\\
    $\{2,4\}$ & present in both source and target\\
    $\{0,3\}$ & absent because the last coordinate of $\deg x^\alpha$ must be positive\\
    $\{1,4\}$ & absent because the first coordinate of $\deg x^\alpha$ must be positive\\
    $\{1,3,4\}$ & absent because the first coordinate of $\deg x^\alpha$ must be positive\\
    $\{0,2,4\}$ & absent because the sum of the coordinates of $\deg x^\alpha$ is at most $-3$\\
    $\{0,1,3\}$ & absent because the last coordinate of $\deg x^\alpha$ must be positive\\
    $\{1,2,4\}$ & absent because the sum of the last 2 coordinates of $\deg x^\alpha$ is at most $-2$\\
    $\{0,2,3\}$ & absent because the sum of the first 2 coordinates of $\deg x^\alpha$ is at most $-2$\\
  \end{tabular}
  \caption{Elements in $H_B^2(S)$ of degrees in the source and target of $\varphi$.}
  \label{fig:cohomology-types}
\end{figure}

According to Figure~\ref{fig:cohomology-types} we can limit our consideration to Laurent monomials $x^\alpha$ with $\Neg\alpha=\{0,2\}$ or $\Neg\alpha=\{2,4\}$.  Since the entries of $\varphi$ are given by multiplication with $x_0$ and $x_4$, the image of $x^\alpha$ cannot acquire additional negative exponents.  Thus we can consider $\varphi$ on each type of cohomology separately, and by symmetry the first is enough:  Laurent monomials $x^\alpha$ with $\Neg\alpha=\{0,2\}$.

We first show that multiplication by $x_4$ is an isomorphism on these monomials in degrees $(-1,0,0),(0,-1,0),(0,0,-1)$.  Note that the coordinates of the degrees of $\frac{1}{x_0},x_1,\frac{1}{x_2},x_3$ all sum to $0$ or $-1$.  The coordinates of the twists appearing in the target of $\varphi$ sum to $-a+m-b=0$, so at least one $x_4$ is required to attain degrees $(-1,0,0),(0,-1,0),(0,0,-1)$.  Thus multiplication by $x_4$ is surjective, and it is injective because $x_4$ does not appear in denominators. Multiplication by $x_4$ pairs the source summand indexed by $(a,b)$ with the target summand indexed by $(a+1,b)$. The unpaired target summand has first index $0$ and vanishes: its first degree coordinate is $0$ or $-1$, whereas $\alpha_0-\alpha_1+\alpha_2\leq-2$ whenever $\Neg\alpha=\{0,2\}$.

Order the Laurent monomial generators of this type in the source and target of $\varphi$ so that $x_4$ acts by the identity matrix and $x_0$ appears with increasingly negative powers.  Then $\varphi$ acts by an upper triangular matrix with 1s along the diagonal, so it is an isomorphism as well.  Thus we conclude \eqref{it:phi-iso} from Claim B.

Finally, since $I^m$ is $B$-saturated, $H_B^1(I^m)=0$. The long exact sequence therefore shows that $\varphi$
is injective in degree $(0,0,0)$, finishing Claim B.

\end{proof}

\appendix

\section{Properties of Painlevé--Kuratowski Convergence}\label{appendix}

In this appendix we provide proofs for the main properties of Painlevé--Kuratowski convergence used at various points throughout. These results are relatively straightforward and stated elsewhere but included here for completeness.

\begin{lemma}\label{lem:unpacking-region-limit-defs}
Let $\{A_{p}\}$ be a sequence of subsets of $\RR^r$.
\begin{enumerate}
\item A point $x\in\lcl \{A_{p}\}$ if and only if for all $\epsilon>0$ there exists a constant $N\in \ZZ_{\geq1}$ such that if $p\geq N$ then $d(x,A_{p})\leq \epsilon$.
\item A point $x\in \ucl \{A_{p}\}$ if and only if for all $\epsilon>0$ there are infinitely many $p\in \ZZ_{\geq1}$ such that $d(x,A_{p})\leq \epsilon$.
\end{enumerate}
\end{lemma}

\begin{proof}
  Both claims follow from the following characterization of the limit superior: if $\{d_{p}\}$ is sequence of real numbers then the limit superior of $\{d_{p}\}$ is equal to $L\in \RR$ if and only if both
  \begin{enumerate}
  \item for all $\epsilon>0$ there exists a constant $N$ such that if $p\geq N$ then $d_{p} \leq L+\epsilon$, and
  \item for all $\epsilon>0$ the set $\{p \in \ZZ_{\geq1} \; | \; d_{p} \geq L-\epsilon\}$ is infinite.
  \end{enumerate}
  In particular, if $L=0$ and $d_{p}\geq 0$ for all $p$ then the second condition, that $d_{p}\geq-\epsilon$ for infinitely many $p$, is vacuous since $d_{p}\geq0>-\epsilon$ for all $p$. Part~(1) of the Lemma is then the first condition for $d_{p}=d(x,A_{p})$.

  The second part of the lemma follows in a similar fashion. Using that $\liminf \{d_{p}\}=-\limsup \{-d_{p}\}$, the limit inferior of $\{d_{p}\}$ is equal to $L\in \RR$ if and only if both
  \begin{enumerate}
  \item for all $\epsilon>0$ there exists a constant $N$ such that if $p\geq N$ then $d_{p} \geq L-\epsilon$, and
  \item for all $\epsilon>0$ the set $\{p \in \ZZ_{\geq1} \; | \; d_{p}\leq L+\epsilon\}$ is infinite.
  \end{enumerate}
  Again the first condition is vacuous and the second gives part~(2) of the Lemma.
\end{proof}

\begin{lemma}\label{lem:region-limit-properties}
Let $\{A_{p}\}$ be a sequence of subsets of $\RR^r$.
\begin{enumerate}
\item The upper closed limit $\ucl \{A_{p}\}$ is a closed subset of $\RR^r$.
\item The lower closed limit $\lcl \{A_{p}\}$ is a closed subset of $\RR^r$.
\item There is an inclusion
\[
\lcl \{A_{p}\} \subseteq \ucl \{A_{p}\}.
\]
\item Upper closed limits and lower closed limits are independent of taking closures:
\[
\ucl \{A_{p}\}=\ucl \{\overline{A}_{p}\} \quad \text{and} \quad \lcl \{A_{p}\}=\lcl \{\overline{A}_{p}\}.
\]
\item If the sequence $\{A_{p}\}$ is constant, i.e., $A_{p}=A$ for all $p$ then $\lim \{A_{p}\} = \overline{A}$.
\item If $\{B_{p}\}$ is a sequence of subsets such that $A_{p}\subseteq B_{p}$ for all $p$ then
\[
\ucl \{A_{p}\} \subseteq \ucl \{B_{p}\}.
\]
\item If $B \subseteq \RR^r$ is a closed subset such that $A_{p} \subseteq B$ for all $p$ then
\[
\lcl \{A_{p}\} \subseteq \ucl \{A_{p}\} \subseteq B.
\]
\item If $C\subseteq \RR^r$ such that $C\subseteq A_{p}$ for all $p$ then
\[
C\subseteq \lcl \{A_{p}\} \subseteq \ucl \{A_{p}\}.
\]
\end{enumerate}
\end{lemma}

\begin{proof}

  For ease of notation throughout this proof we set $L=\lcl \{A_{p}\}$ and $U=\ucl \{A_{p}\}$.

  We begin with part~(1). Let $x\in \overline{U}$.  By the definition of upper closed limits $x\in U$ if and only if $\liminf d(x,A_{p})=0$. This is equivalent to showing that for all $\epsilon>0$ there are infinitely many $p$ such that $d(x,A_{p})\leq \epsilon$. As $x$ is in the closure of $U$ there exists $y\in U$ such that $d(x,y)\leq \frac{\epsilon}{2}$. Thus by the triangle inequality we have that
  \[
  d(x,A_{p})\leq d(x,y)+d(y,A_{p})\leq \frac{\epsilon}{2}+d(y,A_{p}).
  \]
  However, by Lemma~\ref{lem:unpacking-region-limit-defs} $y\in U$ if and only if $d(y,A_{p})\leq \frac{\epsilon}{2}$ for infinitely many $p$. Thus the above inequality shows that $d(x,A_{p})\leq \epsilon$ for infinitely many $p$, so $x\in U$.

  Part (2) is similar to part~(1). Let $x\in \overline{L}$. We show that $x\in L$ by proving that for all $\epsilon>0$ there exists $N$ such that $d(x,A_{p})\leq \epsilon$ for all $p\geq N$. Fix $\epsilon>0$. Since $x$ is in the closure of $L$ there exists $y\in L$ such that $d(x,y)\leq \frac{\epsilon}{2}$. As $y\in L$ we know that there exists a constant $N$ such that $d(y,A_{p})\leq \frac{\epsilon}{2}$ for all $p\geq N$. The result now follows from the triangle inequality, as for all $p\geq N$
  \[
  d(x,A_{p})\leq d(x,y)+d(y,A_{p})\leq \epsilon.
  \]

  Part (3) is immediate from Lemma~\ref{lem:unpacking-region-limit-defs} since if for all $\epsilon>0$ there exists a constant $N$ such that $d(x,A_{p})\leq \epsilon$ (i.e., $x\in L$) then $d(x,A_{p})\leq \epsilon$ for infinitely many $p$ (i.e., $x\in U$).

  Turning to part~(4), the inclusions $\ucl\{A_p\}\subseteq\ucl\{\overline A_p\}$ and $\lcl\{A_p\}\subseteq\lcl\{\overline A_p\}$ hold by definition.  Let $x\in\lcl\{\overline A_p\}$ and fix $\epsilon>0$.  Then by Lemma~\ref{lem:unpacking-region-limit-defs} there exists $N\in\ZZ_{\geq 1}$ such that $d(x,\overline A_p)\leq\frac{\epsilon}{2}$ for $p\geq N$.  Choose a sequence $y_N,y_{N+1},\ldots$ with $y_p\in\overline A_p$ witnessing $d(x,y_p)\leq\frac{\epsilon}{2}$.  For each $p$ there exists $y'_p\in A_p$ with $d(y_p,y_p')<\frac{\epsilon}{2}$ by definition of closure.  Thus
  \[d(x,A_p)\leq d(x,y_p)\leq d(x,y_p)+d(y_p,y'_p)\leq\epsilon\]
  by the triangle inequality so $x\in\lcl\{A_p\}$.  The same argument applied to a subsequence of $\{A_p\}$ shows $\ucl\{\overline A_p\}\subseteq\ucl\{A_p\}$.

  We now prove part~(5). By part~(4), without loss of generality we may assume that $A$ is closed. Since the sequence $\{A_p\}$ is constant we see immediately from Lemma~\ref{lem:unpacking-region-limit-defs} that $\lcl \{A_{p}\}=\ucl \{A_{p}\}$, and $x\in \lim \{A_{p}\}$ if and only if $d(x,A)\leq \epsilon'$ for all $\epsilon'>0$.  The set of $x$ satisfying this condition is exactly $\overline A$.

  We now prove part~(6). Let $x\in \ucl\{A_{p}\}$. By definition this means $\liminf d(x,A_{p})=0$. Since $A_{p}\subseteq B_{p}$ we have $d(x,B_{p})\leq d(x,A_{p})$ for all $p$, so by a property of the limit inferior
  \[
  \liminf\limits_{p\to\infty}d(x,B_{p})\leq \liminf\limits_{p\to\infty}d(x,A_{p})=0.
  \]
  However, since $d(x,B_{p})\geq0$ for all $p$ this means that $\liminf d(x,B_{p})=0$, which implies that $x\in \ucl \{B_{p}\}$.

  Part (7) follows immediately from parts (5) and (6) by letting $B_{p}=B$ be a constant sequence of subsets.

  Finally for part~(8), if $C\subseteq A_{p}$ for all $p$ then for all $x\in C$ we have $d(x,A_{p})=d(x,x)=0$ for all $p$. Hence $\limsup d(x,A_{p})=0$, implying $x\in \lcl\{A_{p}\}$ as needed.
\end{proof}

\begin{lemma}\label{lem:region-chains}
Let $\{A_{p}\}$ be a sequence of subsets of $\RR^r$.
\begin{enumerate}
\item If $A_{1}\subseteq A_{2}\subseteq A_{3}\subseteq \cdots$ then $\lim\{A_{p}\}=\overline{\bigcup A_{p}}$.
\item If $A_{1}\supseteq A_{2}\supseteq A_{3}\supseteq \cdots$ then $\lim\{A_{p}\}=\bigcap \overline{A}_p$
\end{enumerate}
\end{lemma}

\begin{proof}[Proof of Lemma~\ref{lem:region-chains}]
  By part~(4) of Lemma \ref{lem:region-limit-properties}, without loss of generality we may assume that $A_{p}$ is closed for all $p$.

  We begin with part~(1). We will show that
  \[\textstyle
  \overline{\bigcup A_p}\subseteq \lcl \{A_{p}\} \subseteq \ucl \{A_{p}\} \subseteq \overline{\bigcup A_p}.
  \]
  The middle inclusion is part~(3) of Lemma~\ref{lem:region-limit-properties}, while the inclusion on the right side is part~(7) of Lemma~\ref{lem:region-limit-properties}.  Fix $p'\geq 1$.  Since the definition of $\lcl\{A_p\}$ only depends on $A_p$ for $p\gg 0$ we may replace $\{A_p\}_{p\geq 1}$ by $\{A_p\}_{p\geq p'}$ without affecting the limit.  Thus $\lcl\{A_p\}$ contains $A_{p'}$ by part~(8) of Lemma~\ref{lem:region-limit-properties}.  Since $\lcl\{A_p\}$ contains $A_{p'}$ for each $p'$ it contains their union and thus the closure by part~(2) of Lemma~\ref{lem:region-limit-properties}.

  For part~(2) we need
  \[\textstyle
  \bigcap\overline{A}_p\subseteq \lcl \{A_{p}\} \subseteq \ucl \{A_{p}\} \subseteq \bigcap\overline{A}_p.\]
  Again the inclusions follow from Lemma~\ref{lem:region-limit-properties}: the middle is part~(3) of and the left is part~(8).  Fixing $p'$ and replacing $\{A_p\}_{p\geq 1}$ with $\{A_p\}_{p\geq p'}$ we see that $\ucl\{A_p\}\subseteq\overline{A}_{p'}$ by part~(7).
\end{proof}

\bibliographystyle{alpha}
\bibliography{references.bib}

\end{document}